\documentclass[12pt]{amsart}
\usepackage{graphicx}
\usepackage{amssymb}
\usepackage{color}
\usepackage[all]{xy}

\newtheorem{thm}{Theorem}
\newtheorem{lem}{Lemma}

\newtheorem{cor}{Corollary}

\newtheorem{defn}{Definition}

\newtheorem{prp}{Proposition}
\theoremstyle{remark}
\newtheorem{remark}[thm]{Remark}
\DeclareMathOperator{\Des}{Des}

\DeclareMathOperator{\Pe}{Pk}
\DeclareMathOperator{\Sh}{Sh}
\DeclareMathOperator{\Sol}{Sol}

\DeclareMathOperator{\Q}{\mathcal{Q}\textit{sym}}

\title[Colored posets and quasisymmetric functions]{Colored posets and colored quasisymmetric functions}

\author[S. K. Hsiao]{Samuel K. Hsiao}

\author[T. K. Petersen]{T. Kyle Petersen}

\begin{document}
\begin{abstract}
The colored quasisymmetric functions, like the classic quasisymmetric functions, are known to form a Hopf algebra with a natural peak subalgebra. We show how these algebras arise as the image of the algebra of colored posets. To effect this approach we introduce colored analogs of $P$-partitions and enriched $P$-partitions. We also frame our results in terms of Aguiar, Bergeron, and Sottile's theory of combinatorial Hopf algebras and its colored analog.
\end{abstract}

\maketitle

\section{Introduction}

Hopf algebras have a natural role in combinatorics because, in essence, they encode how to assemble and disassemble discrete structures. Examples include the Malvenuto-Reutenauer algebra of permutations, the Loday-Ronco algebra of planar binary trees, and the classic symmetric functions. The interest in Hopf algebras is both reflected in and heightened by the work of Marcelo Aguiar, Nantel Bergeron, and Frank Sottile \cite{AguiarBergeronSottile}, who define the category of \emph{combinatorial Hopf algebras}. The algebra of quasisymmetric functions, $\Q$, is the terminal object in this category. Within $\Q$ lies $\mathbf{\Pi}$, the algebra of peak functions, which is the terminal object in the category of \emph{odd} combinatorial Hopf algebras. Recently, Bergeron and Christophe Hohlweg \cite{BergeronHohlweg} suggested ``colored" versions of combinatorial Hopf algebras.

In this work we take the viewpoint implicit in Claudia Malvenuto's thesis \cite[4.2-4.3]{Malvenuto} that, while perhaps not fundamental in a categorical sense, the algebra of partially ordered sets is nonetheless key to understanding the Hopf structure on $\Q$ and $\mathbf{\Pi}$. The method for moving from posets to quasisymmetric functions lies with the theory of $P$-partitions: combinatorial objects that carry important algebraic structure with them. Here we use the term ``$P$-partition" to refer to Richard Stanley's tool \cite{Stanley1} as well as John Stembridge's \emph{enriched} $P$-partitions \cite{Stembridge}.

In Section \ref{sec:type A} we present the Hopf algebra of labeled posets, $\mathcal{P}$, followed by a quick overview of $\Q$. We then survey Stanley's $P$-partitions and show how simply derived combinatorial lemmas show that $\Q$ is a homomorphic image of $\mathcal{P}$. We follow this result with an analogous treatment of the peak functions. Here Stembridge's enriched $P$-partitions provide the link between $\mathcal{P}$ and $\mathbf{\Pi}$. We conclude the first section by connecting what we have done to the theory of combinatorial Hopf algebras.

Section \ref{sec:color} begins with the definition of the Hopf algebra of \emph{colored} posets, followed by a discussion of colored quasisymmetric functions and peak functions that parallels Section \ref{sec:type A}. Along the way we outline a theory of colored $P$-partitions and colored enriched $P$-partitions, concluding with links to colored combinatorial Hopf algebras. Although we follow the spirit of Bergeron and Hohlweg's work in this regard \cite[Section~5.3]{BergeronHohlweg}, there is some ambiguity in their definition. See Remark \ref{rem:err}. We clarify their definition and present its consequences in Section \ref{sec:colorHopf}.

Throughout the paper, by a \emph{Hopf algebra}, $\mathcal{H}$, we mean a graded connected Hopf algebra $\mathcal{H} = \bigoplus_{n \geq 0} \mathcal{H}_n$, where being connected means $\mathcal{H}_0$ is the one-dimensional vector space spanned by the unit element of $\mathcal{H}$. Further, we assume for simplicity that our algebras are over $\mathbb{Q}$, though many constructions will work over the integers. The assumption that we work over the field of rational numbers makes it simpler to check whether something is indeed a Hopf algebra.

We mention that this work originally arose in an attempt to outline the Hopf structure of the type B quasisymmetric functions and the type B peak functions \cite{Chow,Petersen2}. Motivated by the type A case, the hope was to define a Hopf algebra of type B posets whose image, via Chak-On Chow's type B $P$-partitions \cite{Chow}, would be the type B quasisymmetric functions, and whose image via the second author's type B enriched $P$-partitions \cite{Petersen2} would be the type B peak functions. While we were able to show that the type B quasisymmetric functions and type B peak functions are Hopf algebras, we failed to find an elegant poset structure beneath. Our results can be found in \cite{HsiaoPetersen}.

We conclude the introduction with Table \ref{table:dimension}, which lists the dimensions of the $n$th graded components of the algebras discussed in the paper, apart from the poset algebras (we would be very interested if someone could give us formulas for those dimensions!).
\begin{table}[h]
\vspace{.5cm}
\begin{tabular}{c|c}
  \textbf{Hopf algebra $\mathcal{H}$} & \textbf{Dimension of $\mathcal{H}_n$} \\
\hline \\
 $\Q$ & $2^{n-1}$ = number of compositions of $n$ \\
 & \\
\hline \\
 $\mathbf{\Pi}$ & $f_n$ given by $f_1 = f_2 = 1$, $f_n = f_{n-1} + f_{n-2}$ for $n > 2$ \\
& \\
\hline \\
 $\Q^{(m)}$ & $m(m+1)^{n-1}$ = number of $m$-colored compositions of $n$ \\
& \\
\hline \\
 $\mathbf{\Pi}^{(m)}$ & $f_{m,n}$ given by $f_{m,1} = m$, $f_{m,2} = m^2$, $f_{m,n} = mf_{m,n-1} + f_{m,n-2}$ for $n > 2$
\end{tabular}
\vspace{.5cm}
\caption{Dimensions of $n$th homogeneous components. \label{table:dimension}}
\end{table}

\section{The classical case}\label{sec:type A}

\subsection{Hopf algebra of posets}

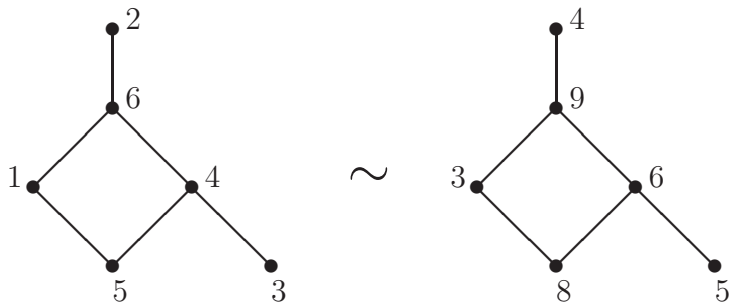
\begin{figure}
\begin{picture}(150,150) \thicklines
   \put(30,60){\line(1,-1){30}} \put(30,60){\line(1,1){30}}
\put(90,60){\line(-1,-1){30}} \put(90,60){\line(-1,1){30}}
\put(90,60){\line(1,-1){30}} \put(60,90){\line(0,1){30}}
\put(30,60){\circle*{5}} \put(60,30){\circle*{5}} \put(60,90){\circle*{5}} \put(60,120){\circle*{5}} \put(90,60){\circle*{5}} \put(120,30){\circle*{5}}
\put(20,60){1} \put(60,17){5} \put(65,90){6} \put(65,120){2} \put(95,60){4} \put(120,17){3}
\end{picture}
\begin{picture}(10,150)
\put(-5,60){\LARGE $\sim$}
\end{picture}
\begin{picture}(150,150) \thicklines
   \put(30,60){\line(1,-1){30}} \put(30,60){\line(1,1){30}}
\put(90,60){\line(-1,-1){30}} \put(90,60){\line(-1,1){30}}
\put(90,60){\line(1,-1){30}} \put(60,90){\line(0,1){30}}
\put(30,60){\circle*{5}} \put(60,30){\circle*{5}} \put(60,90){\circle*{5}} \put(60,120){\circle*{5}} \put(90,60){\circle*{5}} \put(120,30){\circle*{5}}
\put(20,60){3} \put(60,17){8} \put(65,90){9} \put(65,120){4} \put(95,60){6} \put(120,17){5}
\end{picture}
\caption{Two equivalently labeled posets.}\label{fig:equivposets}
\end{figure}

A \emph{labeled poset} is a finite set $P$ with partial order $<_P$ whose elements are distinct positive integers. Two labeled posets $P$ and $Q$ are equivalent, written $P\sim Q$, if there is an isomorphism of posets $\phi: P \to Q$ such that for all $i <_P j$, we have $\phi(i) <_{\mathbb{P}} \phi(j)$ if and only if $i <_{\mathbb{P}} j$ (where $<_{\mathbb{P}}$ is the ordering on $\mathbb{P}$, the set of positive integers). See Figure \ref{fig:equivposets}. Let $\mathcal{P}_n$, for $n\ge 0$, denote the vector space over $\mathbb{Q}$ with basis consisting of all labeled posets of cardinality $n$, modulo the equivalence relation. We will describe a Hopf algebra structure on the graded vector space $\mathcal{P} = \bigoplus_{n \geq 0} \mathcal{P}_n$.

Given two labeled posets $P$ and $Q$, let $P \sqcup Q$ denote their disjoint union as posets. If $P$ and $Q$ have any elements in common, we can always replace $Q$ with a label-equivalent poset that has no elements in common with $P$. For example, since these are finite sets, we can suppose the elements of $Q$ are sufficiently large. It is easy to see that $P \sqcup Q \sim Q \sqcup P$ and if both $P\sim P'$ and $Q\sim Q'$, then $P \sqcup Q\sim P' \sqcup Q'$. Therefore $\sqcup$ is a well-defined commutative product on $\mathcal{P}$. With this product $\mathcal{P}$ becomes a graded algebra whose unit element is the empty poset $\emptyset$.

To define the coproduct, first recall that an \emph{order ideal} (or \emph{lower ideal}) of a poset $P$ is a subset $I$ of $P$ such that if $x\in I$ and $y <_P x$ then $y \in I$. Let $\mathcal{I}(P)$ be the set of order ideals of $P$. For example,
\[ \mathcal{I}\left(
\begin{picture}(50,45)(-5,10)
\put(0,0){\line(1,1){20}} \put(20,20){\line(0,1){20}} \put(20,20){\line(1,-1){20}}
\put(0,0){\circle*{4}} \put(20,20){\circle*{4}} \put(20,40){\circle*{4}} \put(40,0){\circle*{4}}
\put(0,-12){1} \put(25,20){4} \put(25,40){2} \put(40,-12){3}
\end{picture}
 \right)
= \left\{ \emptyset  ,
\begin{picture}(15,5)(-5,-5)
\put(2,0){\circle*{4}} \put(0,-12){1}
\end{picture}
,
\begin{picture}(15,5)(-5,-5)
\put(2,0){\circle*{4}} \put(0,-12){3}
\end{picture}
,
\begin{picture}(35,5)(-5,-5)
\put(2,0){\circle*{4}} \put(0,-12){1} \put(22,0){\circle*{4}} \put(20,-12){3}
\end{picture}
,
\begin{picture}(55,25)(-5,5)
\put(2,0){\circle*{4}} \put(0,-12){1} \put(22,20){\circle*{4}} \put(25,20){4} \put(42,0){\circle*{4}} \put(40,-12){3} \put(2,0){\line(1,1){20}} \put(22,20){\line(1,-1){20}}
\end{picture}
,
\begin{picture}(55,45)(-5,10)
\put(2,0){\circle*{4}} \put(0,-12){1} \put(22,20){\circle*{4}} \put(25,20){4} \put(42,0){\circle*{4}} \put(40,-12){3} \put(22,40){\circle*{4}} \put(25,40){2} \put(2,0){\line(1,1){20}} \put(22,20){\line(1,-1){20}} \put(22,20){\line(0,1){20}}
\end{picture}
\right\}.\]
Every order ideal $I$ and its complement $P\setminus I$ are thought of as labeled subposets of $P$. The coproduct $\delta: \mathcal{P} \to \mathcal{P} \otimes \mathcal{P}$ is defined by
\[ \delta(P) = \sum_{I \in \mathcal{I}(P)} I \otimes (P\setminus I)\]
and the counit is projection onto $\mathcal{P}_0 = \mathbb{Q}$.

The product $\sqcup$ and coproduct $\delta$, along with the unit and counit, give $\mathcal{P}$ the structure of a graded connected Hopf algebra. This algebra was presented in \cite[Example~2.3]{AguiarBergeronSottile}. We should be clear that $\mathcal{P}$ is distinct from Rota's Hopf algebra of (isomorphism classes of finite) graded posets \cite{Rota, JoniRota, Ehrenborg}. 

An explicit formula for the antipode of $\mathcal{P}$ can be deduced from the following general inductive formula due to Milnor and Moore \cite[Section~8]{MilnorMoore}. (Note that they use the term ``conjugate" instead of antipode. See also \cite[Lemma~2.1]{Ehrenborg}.) For any Hopf algebra $\mathcal{H} = \bigoplus_{n \ge 0} \mathcal{H}_n$, the antipode $S:\mathcal{H} \to \mathcal{H}$ is determined by the inductive formula $S(1_{\mathcal{H}}) = 1_{\mathcal{H}}$ and
\begin{equation}\label{eq:antipode}
S(h) = -h - \sum S(h_1) h_2 \quad \text{ for $h \in \mathcal{H}_n, n > 0$,}
\end{equation}
where the coproduct of $h$ is $1\otimes h + h \otimes 1 + \sum h_1 \otimes h_2$. Using this formula it is straightforward to show that the antipode of $\mathcal{P}$ satisfies
\begin{equation}\label{eq:Pantipode}
S(P) = \sum_{{k \ge 1 \atop \emptyset = I_0 \subsetneq I_1 \subsetneq \cdots \subsetneq I_k = P} \atop  I_i \in \mathcal{I}(P) } (-1)^k (I_1\setminus I_0) \sqcup \cdots \sqcup (I_k \setminus I_{k-1})
\end{equation}
for any non-empty labeled poset $P$.

We remark that any permutation $\pi \in \mathfrak{S}_n$ can be considered
a poset with the total order $\pi_i <_{\pi} \pi_{i+1}$, $i=1,2,\ldots,n-1$. For
example, the permutation $\pi = 3214$ is the totally ordered chain $3 <_{\pi} 2
<_{\pi} 1 <_{\pi} 4$. For any poset $P$ of size $n$, let $\mathcal{L}(P)$ denote the Jordan H\"{o}lder set: the set of all totally ordered chains on the elements of $P$ whose relations are consistent with those of $P$. This set is sometimes called the set of ``linear extensions" of $P$. For example let $P$ be the poset defined by $1 >_{P} 5 <_{P} 4$. In this case, 1 and 4 are not comparable, so we have exactly two ways of linearizing $P$: $5 < 4 < 1$ and $5 < 1 < 4$. These correspond to the permutations $541$ and $514$. In pictures,
\[ \mathcal{L}\left(
\begin{picture}(55,25)(-5,10)
\put(0,20){\line(1,-1){20}} \put(20,0){\line(1,1){20}}
\put(0,20){\circle*{4}} \put(20,0){\circle*{4}} \put(40,20){\circle*{4}}
\put(-8,20){1} \put(20,-12){5} \put(45,20){4}
\end{picture}
 \right)
= \left\{
\begin{picture}(20,35)(-5,15)
\put(0,0){\circle*{4}} \put(5,0){5} \put(0,20){\circle*{4}} \put(5,20){4} \put(0,40){\circle*{4}} \put(5,40){1} \put(0,0){\line(0,1){40}}
\end{picture}
,
\begin{picture}(20,35)(-5,15)
\put(0,0){\circle*{4}} \put(5,0){5} \put(0,20){\circle*{4}} \put(5,20){1} \put(0,40){\circle*{4}} \put(5,40){4} \put(0,0){\line(0,1){40}}
\end{picture}
\right\}.\]

\begin{remark}
In fact permutations themselves can be seen to be a Hopf algebra. See work of Malvenuto and Reutenauer \cite{MalvenutoReutenauer}; also \cite{AguiarSottile}. The product and coproduct for the Malvenuto-Reutenauer algebra are not totally unrelated to the product and coproduct on $\mathcal{P}$, but, for instance, their product is not commutative.
\end{remark}

\subsection{Quasisymmetric functions}
The ring of quasisymmetric functions is well-known (see \cite[ch. 7.19]{Stanley2}). Recall that a quasisymmetric function is a formal series \[Q(x_1, x_2, \ldots ) \in \mathbb{Q}[[x_1, x_2,\ldots ]] \] of bounded degree such that for any composition $\alpha = (\alpha_1, \alpha_2, \ldots, \alpha_k)$, the coefficient of $x_{1}^{\alpha_1} x_{2}^{\alpha_2} \cdots x_{k}^{\alpha_k}$ is the same as the coefficient of $x_{i_1}^{\alpha_1} x_{i_2}^{\alpha_2} \cdots x_{i_k}^{\alpha_k}$ for all $i_1 < i_2 < \cdots < i_k$. Recall that a composition of $n$, written $\alpha \models n$, is an ordered tuple of positive integers $\alpha = (\alpha_1, \alpha_2, \ldots, \alpha_k)$ such that $|\alpha| = \alpha_1 + \alpha_2 + \cdots + \alpha_k = n$. In this case we say that $\alpha$ has $k$ parts, or $l(\alpha) = k$. We can put a partial order on the set of all compositions of $n$ by refinement. The covering relations are of the form \[ (\alpha_1, \ldots, \alpha_i + \alpha_{i+1}, \ldots, \alpha_k )
  < (\alpha_1, \ldots, \alpha_i, \alpha_{i+1}, \ldots, \alpha_k).\] Let $\Q_n$ denote the set of all quasisymmetric functions homogeneous of degree $n$. Then $\Q := \bigoplus_{n \geq 0} \Q_n$ denotes the graded ring of all quasisymmetric functions, where $\Q_0 = \mathbb{Q}$.

The most obvious basis for $\Q_n$ is the set of \emph{monomial} quasisymmetric functions, defined for any composition $\alpha = (\alpha_1, \alpha_2, \ldots, \alpha_k) \models n$,
\[ M_{\alpha} := \sum_{i_1 < i_2 < \cdots < i_k} x_{i_1}^{\alpha_1} x_{i_2}^{\alpha_2} \cdots x_{i_k}^{\alpha_k}.\]
There are $2^{n-1}$ compositions of $n$, and hence, the graded component $\Q_n$ has dimension $2^{n-1}$ as a vector space. We can form another natural basis with the \emph{fundamental} quasisymmetric functions, also indexed by compositions,
\[ F_{\alpha} := \sum_{ \alpha \leq \beta } M_{\beta},\] since, by inclusion-exclusion we can express the $M_{\alpha}$ in terms of the $F_{\alpha}$:
\[ M_{\alpha} = \sum_{ \alpha \leq \beta } (-1)^{l(\beta) - l(\alpha)} F_{\beta}.\]
Alternatively, we can write \[F_{\alpha} = \sum_{ i_1 \leq i_2 \leq \cdots \leq i_n} x_{i_1} x_{i_2} \cdots x_{i_n}, \] where if $s = \alpha_1 + \cdots + \alpha_r$, for some $r=1,\ldots,l(\alpha)$ then $i_s < i_{s+1}$.
As an example, \[ F_{21} = M_{21} + M_{111} = \sum_{i < j} x_i^2 x_j + \sum_{i < j < k} x_i x_j x_k  = \sum_{i \leq j < k} x_i x_j x_k .\] When $n = 0$ we have $M_{()} = F_{()} = 1$.

It has been shown that $\Q$ is a Hopf algebra \cite{Ehrenborg, MalvenutoReutenauer} with the usual product of formal power series, counit that takes functions to their constant term, and the following coproduct, \[\Delta : \Q \to \Q \otimes \Q,\] given first in terms of the monomial basis:
\[ \Delta( M_{\alpha} ) = \sum_{ \beta \gamma = \alpha} M_{\beta} \otimes M_{\gamma},\] where $\beta \gamma = ( \beta_1, \ldots, \beta_{l(\beta)}, \gamma_1, \ldots, \gamma_{l(\gamma)})$ is the \emph{concatenation} of $\beta$ and $\gamma$.

The coproduct is easy to understand in terms of the following operation. Let $\mathbb{P} + \mathbb{P}'$ denote the set $\{1,2,\ldots\} \cup \{1',2',\ldots\}$ with total order \[ 1 < 2 < \cdots < 1' < 2' < \cdots,\] and let $X+Y$ define the set of commuting variables $\{ x_s : s \in \mathbb{P}+\mathbb{P}'\}$, with the convention $x_{i'} = y_i$. For any quasisymmetric function $Q$, the coproduct is equivalent to the map $Q(X) \mapsto Q(X+Y)$. To be precise, we can write $Q(X+Y) = \sum R(X)S(Y)$ for some pairs of quasisymmetric functions $R$ and $S$. Then we define $\Delta(Q) = \sum R \otimes S$. As an example,
\begin{align*}
M_{21}(X+Y) & = \sum_{ s < t \in \mathbb{P}+\mathbb{P}'} x^2_s x_t = \sum_{ i < j} x_i^2 x_j + \sum_{i,j} x_i^2 y_j + \sum_{i < j} y_i^2 y_j \\
 & = M_{21}(X) + M_{2}(X)M_1(Y) + M_{21}(Y),
\end{align*}
and thus \[ \Delta(M_{21}) = M_{21} \otimes 1 + M_2 \otimes M_1 + 1 \otimes M_{21}.\]

There is the following nice formula for the antipode \cite{Ehrenborg, MalvenutoReutenauer}:
\begin{align}\label{eq:antipodeM}
S: \Q & \to \Q, \nonumber \\
M_{\alpha} & \mapsto (-1)^{l(\alpha)}\sum_{\beta \leq \alpha} M_{\overleftarrow{\beta}},
\end{align}
where $\overleftarrow{\beta} = (\beta_k, \beta_{k-1}, \ldots, \beta_1)$, the composition formed by writing the parts of $\beta$ backwards.

\subsection{$P$-partitions}\label{sec:ppartitions}

One interesting application for the quasisymmetric functions is as a way of studying descent classes of permutations. A link between these two ideas is given by a combinatorial tool called $P$-partitions. While their definition is due to Stanley \cite{Stanley1}, it was Ira Gessel who used them to connect permutation descents and quasisymmetric functions \cite{Gessel}.

Compositions encode descent sets in the following way. Recall that a \emph{descent} of a permutation $\pi$ of $[n] := \{ 1,2,\ldots,n\}$ is a position $i \in [n-1]$ such that $\pi_i > \pi_{i+1}$, and that an \emph{increasing run} of a permutation $\pi$ is a maximal subword of consecutive letters $\pi_{i+1} \pi_{i+2} \cdots \pi_{i+r}$ such that $\pi_{i+1} < \pi_{i+2} < \cdots < \pi_{i+r}$. By maximality, we have that if $\pi_{i+1} \pi_{i+2} \cdots \pi_{i+r}$ is an increasing run, then $i$ is a descent of $\pi$ (if $i\neq 0$), and $i+r$ is a descent of $\pi$ (if $i+r \neq n$). For any permutation $\pi$ define the \emph{descent composition}, $C(\pi)$, to be the ordered tuple listing the lengths of the increasing runs of $\pi$. If $C(\pi) = (\alpha_1, \alpha_2, \ldots, \alpha_k)$, we can recover the descent set of $\pi$:
\[ \Des(\pi) = \{ \alpha_1, \alpha_1 + \alpha_2, \ldots, \alpha_1 + \alpha_2 + \cdots + \alpha_{k-1} \}.\] Since $C(\pi)$ and $\Des(\pi)$ have the same information, we can use them interchangeably. For example the permutation $\pi = 345\,26\,1$ (the gaps are intentional) has $C(\pi) = (3,2,1)$ and $\Des(\pi) = \{ 3, 5\}$.

Recall \cite[Chapter 4.5]{Stanley1} that a $P$-partition is an order-preserving map from a poset $P$ to some (countable) totally ordered set. To be precise, let $P$ be any labeled partially ordered set and let $S$ be any totally ordered countable set.

\begin{defn}
The map $f: P \to S$ is a \emph{$P$-partition} if it satisfies the following conditions:
\begin{enumerate}

\item $f(i) \leq f(j)$ if $i <_P j$

\item $f(i) < f(j)$ if $i <_P j$ and $i >_{\mathbb{P}} j$.

\end{enumerate}
\end{defn}

We let $\mathcal{A}(P)$ (or $\mathcal{A}(P;S)$ if we wish to emphasize our choice of image set) denote the set of all $P$-partitions, and encode this set in the generating function \[ \Gamma(P) := \sum_{ f \in \mathcal{A}(P)} \prod_{i\in P} x_{f(i)},\] where $n$ is the number of elements in $P$. If we take $S$ to be the set of positive integers, then it should be clear that $\Gamma(P)$ is always going to be a quasisymmetric function homogeneous of degree $n$. As an easy example, let $P$ be the poset defined by $3 >_P 2 <_P 1$. In this case we have \[ \Gamma(P) = \sum_{ f(3) \geq f(2) < f(1)} x_{f(1)} x_{f(2)} x_{f(3)}.\] Notice that if $P \sim Q$, then $\mathcal{A}(P) = \mathcal{A}(Q)$ and consequently $\Gamma(P) = \Gamma(Q)$.

For chains, i.e., permutations, we have
\begin{align*}
 \mathcal{A}(\pi) = \{ f: \pi \to S  \mid & \,\, f(\pi_1) \leq f(\pi_2) \leq \cdots \leq f(\pi_n) \\
 & \mbox{ and } k \in \Des(\pi) \Rightarrow f(\pi_k) < f(\pi_{k+1})\},
\end{align*}
and
\begin{align*} \Gamma(\pi) &= \sum_{ \substack{ i_1 \leq i_2 \leq \cdots \leq i_n \\ k \in \Des(\pi) \Rightarrow i_k < i_{k+1}}} x_{i_1} x_{i_2} \cdots x_{i_n} \\
 &= F_{C(\pi)},
\end{align*}
so the distinct generating functions for the $\pi$-partitions of permutations $\pi \in \mathfrak{S}_n$ form a basis for $\mathcal{Q}sym_n$. In fact, this is the way Gessel \cite{Gessel} originally defined the quasisymmetric functions.

\begin{remark}
Taking this interpretation of the fundamental quasisymmetric functions---as generating functions for descent classes of permutations---leads to some interesting connections between $\Q_n$ and Solomon's descent algebra $\Sol(A_{n-1})$ \cite{Solomon} (the subalgebra of the group algebra of the symmetric group $\mathfrak{S}_n$ formed by the linear span of sums of permutations with the same descent set). In fact, using Stanley's theory of $P$-partitions, one can show that the two are dual with respect to a certain ``internal" coproduct on $\Q_n$. See \cite{Gessel} and \cite{Petersen1} for more. In this paper we shall be concerned more with the Hopf algebraic structure of $\Q$ and its type B analog. One can carry the duality between $\Q$ and $\Sol := \bigoplus_{n \geq 0} \Sol(A_n)$ to the level of Hopf algebras, see \cite{MalvenutoReutenauer}, but we will not do so here.
\end{remark}

We will list some useful facts from the theory of $P$-partitions, all of which can be found, sometimes implicitly, in \cite{Gessel} or \cite{Stanley1}.

\begin{lem}
The set of all $P$-partitions is the disjoint union of the sets of $\pi$-partitions of all linear extensions $\pi$ of $P$: \[\mathcal{A}(P) = \coprod_{\pi \in \mathcal{L}(P)}
\mathcal{A}(\pi).\]
\end{lem}

\begin{cor}
For any labeled poset $P$ we have
\[\Gamma(P) = \sum_{\pi \in \mathcal{L}(P)} \Gamma(\pi).\]
\end{cor}

\begin{lem}
The set of $(P \sqcup Q)$-partitions is in bijection with the cartesian product of the set of $P$-partitions and the set of $Q$-partitions: \[ \mathcal{A}(P\sqcup Q) \leftrightarrow \mathcal{A}(P) \times \mathcal{A}(Q).\]
\end{lem}

\begin{cor}
For any labeled posets $P$ and $Q$ we have
\begin{equation}\label{eq:Gammamult}
\Gamma(P \sqcup Q) = \Gamma(P)\Gamma(Q).
\end{equation}
\end{cor}

This last fact leads to a well-known formula for the multiplication of quasisymmetric functions in terms of the fundamental basis. If $\sigma$ is a permutation of $[n]$ and $\tau$ is a permutation of $[n+1,n+m]:=\{n+1,n+2,\ldots,n+m\}$, then taking $P = \sigma$, $Q=\tau$ in \eqref{eq:Gammamult} gives
\begin{align}
F_{C(\sigma)} F_{C(\tau)} &= \Gamma(\sigma \sqcup \tau) \nonumber \\
 &= \sum_{ \pi \in \Sh(\sigma,\tau)} F_{C(\pi)},\label{eq:Fmult}
\end{align}
where the sum is over the set $\Sh(\sigma,\tau)$ of all \emph{shuffles} of $\sigma$ and $\tau$, i.e., those permutations of $[n+m]$ whose restriction to the letters $[n]$ is $\sigma$ and whose restriction to $[n+1,n+m]$ is $\tau$. This formula amounts to the fact that the set of linear extensions of the disjoint union of two chains is the set of shuffles of those chains.

\subsection{Viewing $\Q$ as the homomorphic image of $\mathcal{P}$}\label{sec:QfromP}

We can use the $P$-partition generating function to define the map $\Gamma: \mathcal{P} \to \Q$ by sending a poset $P$ to $\Gamma(P)$. Malvenuto showed that $\Gamma$ is a homomorphism of Hopf algebras \cite{Malvenuto}. From the perspective of $P$-partitions this is readily apparent. We already observed that $\Gamma(P) = \Gamma(Q)$ whenever $P\sim Q$, so we may think of $\Gamma$ as a well-defined linear map from $\mathcal{P}$ to $\Q$. Formula \eqref{eq:Gammamult} shows that $\Gamma$ is an algebra homomorphism. To see that $\Gamma$ is a coalgebra homomorphism requires only slightly more thought.

Recall that $\mathbb{P} + \mathbb{P}'$ denotes the set $\{ 1,2,\ldots \} \cup \{1',2',\ldots \}$ with total order \[ 1 < 2 < \cdots < 1' < 2' < \cdots.\]

\begin{lem}\label{lem:coprod}
For any labeled poset $P$ we have a bijection \[ \mathcal{A}(P; \mathbb{P} + \mathbb{P}') \leftrightarrow \coprod_{I \in \mathcal{I}(P)} \mathcal{A}(I; \mathbb{P}) \times \mathcal{A}(P\setminus I; \mathbb{P}').\]
\end{lem}

\begin{proof}
For a given $f \in \mathcal{A}(P; \mathbb{P} + \mathbb{P}')$, let $I = \{ x \in P : f(x) \in \mathbb{P} \}$ and let $f|_I$ denote the restriction of $f$ to $I$ and let $f|_{P\setminus I}$ denote the restriction of $f$ to the complement of $I$. If $x \in I$ and $y <_P x$, then by definition of a $P$-partition, $f(y) \leq f(x)$. This means $f(y) \in \mathbb{P}$ and so $y \in I$. Hence $I$ is an order ideal and it is clear that the map $f \mapsto (f|_I , f|_{P\setminus I})$ is a bijection.
\end{proof}

\begin{cor}
For any labeled poset $P$ we have
\begin{equation}\label{eq:gammacoprod}
\Gamma(P)(X+Y) = \sum_{I \in \mathcal{I}(P)} \Gamma(I)(X)\Gamma(P\setminus I)(Y).
\end{equation}
\end{cor}

Equation \eqref{eq:gammacoprod} gives \[\Delta(\Gamma(P)) = \sum_{I \in \mathcal{I}(P)} \Gamma(I) \otimes \Gamma(P\setminus I) = \Gamma\otimes \Gamma( \delta(P)), \] which shows that $\Gamma$ is a coalgebra homomorphism, and hence is a Hopf homomorphism as claimed.

\begin{thm}
The map $\Gamma$ is a homomorphism of Hopf algebras. Moreover, \[\Gamma(\mathcal{P}) = \Q. \]
\end{thm}

We also get the following nice formula for the coproduct in terms of the fundamental basis by thinking of the order ideals of a chain. For a permutation $\pi = \pi_1 \cdots \pi_n$,
\begin{equation}\label{eq:Fcoproduct}
 \Delta(F_{C(\pi)}) = \sum_{i = 0}^n F_{C(\pi_1 \cdots \pi_i)} \otimes F_{C(\pi_{i+1} \cdots \pi_n)}.
 \end{equation}
There is also a simple formula for the antipode found in \cite{Ehrenborg, MalvenutoReutenauer}. It is easily derived from the inductive formula \eqref{eq:antipode} for the antipode. For any permutation $\pi = \pi_1 \cdots \pi_n$,
\begin{equation}\label{eq:Fantipode}
S(F_{C(\pi)}) = (-1)^n F_{C(\pi_n \pi_{n-1} \cdots \pi_1)}.
\end{equation}
Unlike the situation for product and coproduct, the antipode formula \eqref{eq:Fantipode} is not an immediate consequence of the antipode formula for a chain in $\mathcal{P}$. In Section \ref{sec:colorantipode} we will give a simple proof of an analogous formula for colored quasisymmetric functions.

\subsection{Peak functions}\label{sec:peakA}

There is a Hopf subalgebra of $\Q$ given by Stembridge's \emph{peak functions} \cite{Stembridge}, defined as follows. First, for any $\alpha \models n$, define the more general functions
\begin{equation}\label{eq:KA-M}
K_{\alpha} := \sum_{ \alpha \leq \beta^* } 2^{l(\beta)} M_{\beta},
\end{equation}
where $\beta^*$ is the refinement of $\beta$ obtained by replacing, for $i > 1$, every part $\beta_i \geq 2$ with the two parts $(1,\beta_i -1)$. For example, if $\beta = (2,3,1,2)$, then $\beta^* = (2,1,2,1,1,1)$.

Consider compositions of the form $\alpha = (\alpha_1, \alpha_2, \ldots, \alpha_k) \models n$ such that $\alpha_i > 1$ unless $i = k$. In other words, all but the last part of the composition must be greater than 1. We call such compositions \emph{peak compositions}. Let $\mathbf{\Pi}_n$ denote the span of the functions $K_{\alpha}$, where $\alpha$ is a peak composition. These $K_{\alpha}$ are the peak functions, which are shown in \cite{Stembridge} to be $F$-positive and linearly independent, forming a basis for $\mathbf{\Pi}_n$ of dimension $f_n$, the Fibonacci number defined by $f_1 = f_2 = 1$ and $f_n = f_{n-1} + f_{n-2}$ for $n > 2$. We define $\mathbf{\Pi} := \bigoplus_{n \geq 0} \mathbf{\Pi}_n$ to be the graded ring of peak functions.

For any composition $\alpha$, let $\widehat{\alpha}$ denote the composition formed by replacing consecutive 1s with their sum added to the next part to the right: \[ ( \ldots, \alpha_{i}, \underbrace{1,1,\ldots, 1}_r, \alpha_{i+r+1}, \ldots ) \mapsto (\ldots, \alpha_i, \alpha_{i+r+1}+r, \ldots ),\] where $\alpha_i, \alpha_{i+r+1} > 1$. For example, if $\alpha = (3,1,1,3,2,1,1,1)$, then $\widehat{\alpha} = (3,5,2,3)$. The map $\alpha \mapsto \widehat{\alpha}$ is thus a surjection from all compositions onto peak compositions. In \cite{Stembridge} Stembridge defines the map $\Theta: \Q \to \mathbf{\Pi}$ given by \[ \Theta( F_{\alpha} ) = K_{\widehat{\alpha}}\] and shows that $\Theta$ is a graded, surjective homomorphism of rings.

In \cite{BergeronMSW} it is shown that the $K_{\alpha}$ (for any $\alpha \models n$) satisfy the following coproduct formula:
\begin{equation}\label{eq:peakcoprod}
\Delta(K_{\alpha}) = 1 \otimes K_{\alpha} + \sum_{ \beta (i+j) \gamma = \alpha} K_{\beta i}  \otimes K_{ j \gamma},
\end{equation}
where $i > 0$, $j \geq 0$ are integers, and if $j = 1$ then we understand $j \gamma = (1+\gamma_1, \gamma_2, \ldots)$. In particular, if $\alpha$ is a peak composition, then so are all the compositions appearing on the right hand side of \eqref{eq:peakcoprod}. Thus, the peak functions are closed under coproduct and $\mathbf{\Pi}$ is a Hopf subalgebra of $\Q$. The map $\Theta$ can now be seen to be a Hopf homomorphism.
We derive this result with Stembridge's theory of \emph{enriched} $P$-partitions as well.

\subsection{Enriched $P$-partitions}

Peak functions get their name because, just as the $F_{\alpha}$ encode descent classes of permutations, the $K_{\alpha}$, where $\alpha$ is a peak composition, encode peak classes of permutations. While one can define peak functions independently of peaks of permutations (as we have just done), Stembridge's definition of the peak functions was as generating functions for \emph{enriched} $P$-partitions, which happen to depend on peaks.

An \emph{interior peak} of a permutation $\pi$ of $[n]$ is any $i \in [2,n-1]$ such that $\pi_{i-1} < \pi_i > \pi_{i+1}$. Let $\Pe(\pi)$ denote the set of all interior peaks of $\pi$. Valid peak sets are precisely those subsets $I$ of $[2,n-1]$ with the property that if $i \in I$, $i-1 \notin I$. Thus we see that peak sets correspond to the peak compositions using the same map to move between descent sets and ordinary compositions. Let $\widehat{C}(\pi)$ denote the composition corresponding to $\Pe(\pi)$. Roughly speaking, it records the distances between peaks. For example, if $\pi = 327\,5418\,6$, then $\Pe(\pi) = \{ 3, 7\}$, and $\widehat{C}(\pi) = (3,4,1)$.

Just as with ordinary $P$-partitions, Enriched $P$-partitions are maps from a poset to some totally ordered set, but we impose some extra conditions. First, the image set takes a certain form. For any set $S = \{s_1, s_2, \ldots \}$ with total order $ s_1 < s_2 < \cdots$, define the set $S^{\pm} = \{ s_1, s_2, \ldots \} \cup \{ -s_1, -s_2, \ldots \}$ with total order \[-s_1 < s_1 < -s_2 < s_2 < \cdots.\] We write $s \leq^+ t$ in $S^{\pm}$ if $s <_{S^{\pm}} t$ or if $s = t \in S$. Similarly, $s \leq^{-} t$ means $s <_{S^{\pm}} t$ or $-s = -t \in S$. If $S = \mathbb{P}$, then we have \[ -1 < 1 < -2 < 2 < \cdots.\]

\begin{defn}
An \emph{enriched $P$-partition} is a map $f: P \to S^{\pm}$ such that
\begin{enumerate}
\item $f(i) \leq^+ f(j)$ if $i <_P j$ and $i <_{\mathbb{P}} j$

\item $f(i) \leq^- f(j)$ if $i <_P j$ and $i >_{\mathbb{P}} j$.
\end{enumerate}
\end{defn}

We let $\mathcal{E}(P)$ (or $\mathcal{E}(P;S^{\pm})$) denote the set of all enriched $P$-partitions. The generating function for enriched $P$-partitions is \[ \Lambda(P) := \sum_{ f \in \mathcal{E}(P)} \prod_{i \in P} x_{|f(i)|}.\] If we take $S = \mathbb{P}$ then $\Lambda(P)$ is a quasisymmetric function homogeneous of degree $n$. Also note that if $P \sim Q$ then $\mathcal{E}(P) = \mathcal{E}(Q)$ and $\Lambda(P) = \Lambda(Q)$.

Though not immediately obvious from the definition, \cite[Proposition~2.2]{Stembridge} gives the following fact: \[ \Lambda(\pi) = K_{\widehat{C}(\pi)}.\]

\begin{remark}
The sums of permutations of $[n]$ with the same set of interior peaks span a subalgebra of $\Sol(A_{n-1})$ called the \emph{peak algebra}, $\mathfrak{P}_n$, originally due to K. Nyman \cite{Nyman}. Just as $P$-partitions can be used to show duality between $\Q_n$ and $\Sol(A_{n-1})$, enriched $P$-partitions are employed in \cite{Petersen2} to show duality between $\mathbf{\Pi}_n$ and $\mathfrak{P}_n$.
\end{remark}

Now we present some useful facts from the theory of enriched $P$-partitions, taken from \cite{Stembridge}.

\begin{lem}
The set of all enriched $P$-partitions is the disjoint union of the sets of enriched $\pi$-partitions of all linear extensions $\pi$ of $P$: \[\mathcal{E}(P) = \coprod_{\pi \in \mathcal{L}(P)}
\mathcal{E}(\pi).\]
\end{lem}

\begin{cor}
For any labeled poset $P$,
\begin{equation}\label{eq:Lambdasum}
\Lambda(P) = \sum_{\pi \in \mathcal{L}(P)} \Lambda(\pi).
\end{equation}
\end{cor}

\begin{lem}
The set of enriched $(P \sqcup Q)$-partitions is in bijection with the cartesian product of the set of enriched $P$-partitions and the set of enriched $Q$-partitions: \[ \mathcal{E}(P\sqcup Q) \leftrightarrow \mathcal{E}(P) \times \mathcal{E}(Q).\]
\end{lem}

\begin{cor}
For any labeled posets $P$ and $Q$,
\begin{equation}\label{eq:Lambdamult}
\Lambda(P \sqcup Q) = \Lambda(P)\Lambda(Q).
\end{equation}
\end{cor}

Taking $P$ and $Q$ to be chains in \eqref{eq:Lambdamult} gives the following formula for multiplication of peak functions analogous to \eqref{eq:Fmult}:
\begin{equation}\label{eq:Kmult}
K_{\widehat{C}(\sigma)} K_{\widehat{C}(\tau)} = \sum_{ \pi \in \Sh(\sigma,\tau) } K_{\widehat{C}(\pi)},
\end{equation}
where again the sum is over all shuffles of $\sigma$ and $\tau$.

\subsection{Viewing $\mathbf{\Pi}$ as the homomorphic image of $\mathcal{P}$}

Just as in Section \ref{sec:QfromP} we define the map $\Lambda: \mathcal{P} \to \Q$ by taking a poset to its enriched generating function. We will show that this map is a Hopf homomorphism and therefore its image, $\mathbf{\Pi}$, is a Hopf algebra. First we note that equation \eqref{eq:Lambdamult} shows that $\Lambda$ is an algebra homomorphism.

To see that it is a coalgebra homomorphism, let $\mathbb{P}^{\pm} + (\mathbb{P}^{\pm})'$ denote the set $\{ -1,1,-2,2,\ldots \} \cup \{-1',1',-2',2',\ldots \}$ with total order \[ -1 < 1 < -2 < 2 < \cdots < -1' < 1' < -2' < 2' < \cdots.\]

\begin{lem}
For any labeled poset $P$ we have \[ \mathcal{E}(P; \mathbb{P}^{\pm} + (\mathbb{P}^{\pm})') \leftrightarrow \coprod_{I \in \mathcal{I}(P)} \mathcal{E}(I; \mathbb{P}^{\pm}) \times \mathcal{E}(P\setminus I; (\mathbb{P}^{\pm})').\]
\end{lem}

The structure of the proof is identical to that of Lemma \ref{lem:coprod}.

\begin{cor}\label{cor:Lambdacoprod}
For any labeled poset $P$ we have
\[
\Lambda(P)(X+Y) = \sum_{I \in \mathcal{I}(P)} \Lambda(I)(X)\Lambda(P\setminus I)(Y).
\]
\end{cor}

So we can conclude that $\Lambda$ is a coalgebra homomorphism, and hence $\Lambda$ is a Hopf homomorphism.

\begin{thm}
The map $\Lambda$ is a homomorphism of Hopf algebras. Moreover, \[\Lambda(\mathcal{P}) = \mathbf{\Pi}. \]
\end{thm}

As a special case of Corollary \ref{cor:Lambdacoprod} we can rewrite formula \eqref{eq:peakcoprod} in terms of the ideals of a chain, or permutation, $\pi = \pi_1\cdots \pi_n$:
\begin{equation}\label{eq:Kcoprod}
\Delta(K_{\widehat{C}(\pi)}) = \sum_{i = 0}^n K_{\widehat{C}(\pi_1 \cdots \pi_i)} \otimes K_{\widehat{C}(\pi_{i+1} \cdots \pi_n)}.
\end{equation}

By taking formulas \eqref{eq:Fmult} and \eqref{eq:Kmult} we can see that $\Theta$ is an algebra homomorphism and formulas \eqref{eq:Fcoproduct} and \eqref{eq:Kcoprod} show that it is also a coalgebra homomorphism. Note that we can write $\Lambda$ as $\Lambda = \Theta \circ \Gamma$, the composition of two Hopf homomorphisms. We get the antipode formula
\[ S(K_{\widehat{C}(\pi)}) = (-1)^nK_{\widehat{C}(\pi_n \pi_{n-1} \cdots \pi_1)}\]
by \eqref{eq:Fantipode} and the fact that $S\circ \Theta = \Theta\circ S$.

\subsection{Connections with combinatorial Hopf algebras}\label{sec:combHopf}

In addition to its interpretation in terms of $P$-partitions, the maps $\Gamma$ and $\Lambda$ also have simple interpretations in the context of Aguiar, Bergeron, and Sottile's theory of combinatorial Hopf algebras  \cite{AguiarBergeronSottile}. A \emph{combinatorial Hopf algebra} is a pair $(\mathcal{H}, \varphi)$, where $\mathcal{H}$ is a graded connected Hopf algebra and $\varphi$ is a character, i.e., an algebra homomorphism $\mathcal{H} \to \mathbb{Q}$. The canonical character for $\Q$ is defined by setting $x_1 = 1$ and $x_i = 0$ for $i>1$. Equivalently, in terms of the monomial and fundamental functions,
\[
 \zeta_{\mathcal{Q}}(M_\alpha) = \zeta_{\mathcal{Q}}(F_\alpha) = \begin{cases} 1 & \text{ if $\alpha = ()$ or $\alpha = (n)$ for some $n > 0$,} \\
0 & \text{ otherwise.}
\end{cases}
\]
A fundamental result of \cite[Theorem~4.1]{AguiarBergeronSottile} is that $(\Q, \zeta_{\mathcal{Q}})$ is the terminal object in the category of combinatorial Hopf algebras, i.e., for every combinatorial Hopf algebra $(\mathcal{H}, \varphi)$ there exists a unique homomorphism of graded Hopf algebras $\widehat{\varphi}: \mathcal{H} \to \Q$ such that $\varphi = \zeta_{\mathcal{Q}} \circ \widehat{\varphi}.$

Consider the character $\zeta_{\mathcal{P}}: \mathcal{P} \to \mathbb{Q}$ defined for any labeled poset $P$ by
\begin{equation}\label{eq:zetaP}
\zeta_{\mathcal{P}}(P) = \begin{cases} 1 & \text{if $P$ is naturally labeled}  \\  0 & \text{otherwise,} \end{cases}
\end{equation}
where being naturally labeled means that $i <_{\mathbb{P}} j$ whenever $i<_P j$. By considering linear extensions, it is easy to verify that \[\zeta_{\mathcal{P}} = \zeta_{\mathcal{Q}} \circ \Gamma\] and so by uniqueness $\Gamma=\widehat{\zeta}_{\mathcal{P}}$.

\begin{remark}
The algebra of posets is presented in \cite[Example~2.3]{AguiarBergeronSottile}, but the authors study it with the character $\zeta: \mathcal{P} \to \mathbb{Q}$ given by $\zeta(P) = 1$ for any poset $P$.
\end{remark}

We now try to understand $\Lambda$ from the perspective of combinatorial Hopf algebras. Recall that the characters of a graded connected Hopf algebra $\mathcal{H}$ form a group under convolution product:
\[\zeta \varphi (h) = \sum\zeta(h_1) \varphi(h_2)\]
where the coproduct sends $h\in \mathcal{H}$ to $\sum h_1 \otimes h_2$. The inverse of a character is obtained by composing with the antipode: $\zeta^{-1} = \zeta \circ S$. If $\zeta$ is a character, let $\bar{\zeta}$ denote the character defined by $\bar{\zeta}(h) = (-1)^n \zeta(h)$ for any $h \in \mathcal{H}_n, n \ge 0$. A combinatorial Hopf algebra $(\mathcal{H}, \zeta)$ is called \emph{odd} if $\zeta$ is an odd character, meaning $\bar{\zeta} = \zeta^{-1}$. Given any character $\zeta$, the character $\nu = \bar{\zeta}^{-1}\zeta$ is always odd.

As a combinatorial Hopf algebra, the algebra of peak functions is the pair $(\mathbf{\Pi}, \zeta_{\mathcal{Q}})$. On the $K$-basis, we have:
\begin{equation}\label{eq:Kzeta}
\zeta_{\mathcal{Q}}(K_{\alpha}) =
\begin{cases}
1 & \mbox{if } \alpha = (), \\
2 & \mbox{if } \alpha = (n) \mbox{ for some } n>0,\\
0 & \mbox{otherwise.}
\end{cases}
\end{equation}
Aguiar, Bergeron, and Sottile \cite[Corollary 6.2]{AguiarBergeronSottile} proved that $(\mathbf{\Pi},\zeta_{\mathcal{Q}})$ is the terminal object for odd combinatorial Hopf algebras. They also showed that the map $\Theta:\Q \to \mathbf{\Pi}$ is the unique Hopf homomorphism induced by the odd character $\nu_{\mathcal{Q}} := \bar{\zeta}_{\mathcal{Q}}^{-1}  {\zeta}_{\mathcal{Q}}$. In other words, $\nu_{\mathcal{Q}} = \zeta_{\mathcal{Q}} \circ \Theta$.

Consider the odd character $\nu_{\mathcal{P}} :=  \bar{\zeta}_{\mathcal{P}}^{-1}  {\zeta}_{\mathcal{P}}$. Since $\zeta_{\mathcal{P}} = \zeta_{\mathcal{Q}} \circ \Gamma$ and $\Gamma$ is a Hopf homomorphism, it follows that $\nu_{\mathcal{P}} = \nu_{\mathcal{Q}} \circ \Gamma$, and so
\begin{equation}\label{eq:nu_P}
 \nu_{\mathcal{P}} = \zeta_{\mathcal{Q}} \circ \Lambda.
 \end{equation}
By uniqueness we have $\Lambda = \widehat{\nu}_{\mathcal{P}}$.

Combining \eqref{eq:nu_P}, \eqref{eq:Kzeta}, and \eqref{eq:Lambdasum}, the following is immediate.
\begin{prp}
For any non-empty labeled poset $P$,
\[ \nu_{\mathcal{P}}(P) = 2 \cdot \#\{\pi\in\mathcal{L}(P)~|~\Pe(\pi) = \emptyset\}.\]
\end{prp}

To make the similarity with $\zeta_{\mathcal{P}}$ clear, notice that we could have written \[\zeta_{\mathcal{P}}(P) = \#\{ \pi \in \mathcal{L}(P)~|~\Des(\pi) = \emptyset\}\] in place of \eqref{eq:zetaP}.

\section{The colored case}\label{sec:color}

\subsection{The Hopf algebra of colored posets}

Let $\omega$ be any primitive $m$th root of unity, and let $\mathcal{C}_m = \{ 1,\omega,\ldots, \omega^{m-1}\}$ be the cyclic group of order $m$ generated by $\omega$. Define the set $\mathbb{P}_m := \mathbb{P} \times \mathcal{C}_m$ with total order \[ 1 < \omega 1 < \cdots < \omega^{m-1} 1 < 2 < \omega 2 < \cdots < \omega^{m-1} 2 < \cdots. \]

We call this set the \emph{$m$-colored natural numbers}. In general, for any countable totally ordered set $S = \{ s_1 < s_2 < \cdots \}$, the set $S_m := S \times \mathcal{C}_m$ is totally ordered as \[s_1 < \omega s_1 < \cdots < \omega^{m-1} s_1 < s_2 < \omega s_2 < \cdots < \omega^{m-1} s_2 < \cdots. \] For colored numbers, define by $\varepsilon(\omega^j s_i) = j$ the \emph{color} of $\omega^j s_i$. However, we sometimes use the term ``color" to refer to the entire coefficient $\omega^j$ as well. The context will make the meaning clear.

\begin{defn}
An $m$-colored poset of $n$ elements, or an $(m,n)$-poset, is a poset $P$ whose elements form a subset of $\mathbb{P}_m$ with distinct absolute values.
\end{defn}

\begin{figure}
\begin{picture}(150,150) \thicklines
   \put(30,60){\line(1,-1){30}} \put(30,60){\line(1,1){30}}
\put(90,60){\line(-1,-1){30}} \put(90,60){\line(-1,1){30}}
\put(90,60){\line(1,-1){30}} \put(60,90){\line(0,1){30}}
\put(30,60){\circle*{5}} \put(60,30){\circle*{5}} \put(60,90){\circle*{5}} \put(60,120){\circle*{5}} \put(90,60){\circle*{5}} \put(120,30){\circle*{5}}
\put(12,60){{\color{blue}$\omega 1$ }} \put(60,17){{\color{blue}$\omega 5$}} \put(65,90){6} \put(65,120){{\color{blue}$\omega 2$}} \put(95,60){{\color{red}$\omega^2 4$}} \put(120,17){3}
\end{picture}
\begin{picture}(10,150)
\put(-5,60){\LARGE $\sim$}
\end{picture}
\begin{picture}(150,150) \thicklines
   \put(30,60){\line(1,-1){30}} \put(30,60){\line(1,1){30}}
\put(90,60){\line(-1,-1){30}} \put(90,60){\line(-1,1){30}}
\put(90,60){\line(1,-1){30}} \put(60,90){\line(0,1){30}}
\put(30,60){\circle*{5}} \put(60,30){\circle*{5}} \put(60,90){\circle*{5}} \put(60,120){\circle*{5}} \put(90,60){\circle*{5}} \put(120,30){\circle*{5}}
\put(12,60){{\color{blue}$\omega 3$}} \put(60,17){{\color{blue}$\omega 8$}} \put(65,90){9} \put(65,120){{\color{blue}$\omega 4$}} \put(95,60){{\color{red}$\omega^2 6$}} \put(120,17){5}
\end{picture}
\caption{Two equivalently labeled colored posets.}\label{fig:equiv}
\end{figure}
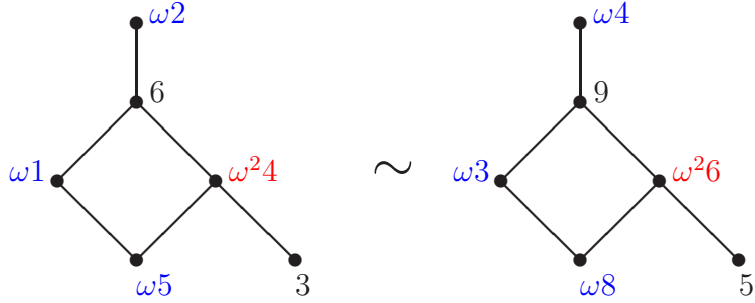

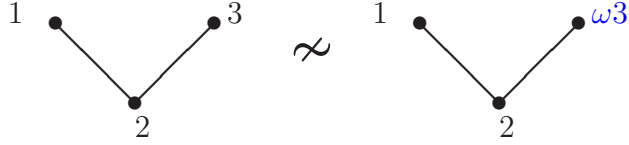
\begin{figure}
\begin{picture}(120,75)(0,10) \thicklines
   \put(30,60){\line(1,-1){30}}
\put(90,60){\line(-1,-1){30}}
\put(30,60){\circle*{5}} \put(60,30){\circle*{5}} \put(90,60){\circle*{5}}
\put(12,60){1} \put(60,17){2}  \put(95,60){3}
\end{picture}
\begin{picture}(10,75)(0,10)
\put(-5,45){\LARGE $\nsim$}
\end{picture}
\begin{picture}(120,75)(0,10) \thicklines
   \put(30,60){\line(1,-1){30}}
\put(90,60){\line(-1,-1){30}}
\put(30,60){\circle*{5}} \put(60,30){\circle*{5}} \put(90,60){\circle*{5}}
\put(12,60){1} \put(60,17){2}  \put(95,60){{\color{blue}$\omega 3$}}
\end{picture}
\caption{Inequivalently labeled colored posets.}\label{fig:inequiv}
\end{figure}

We say that two colored posets $P$ and $Q$ have equivalent labelings, written $P \sim Q$, if there is an isomorphism of posets $\phi:P \to Q$ such that:
\begin{enumerate}
\item the map $\phi$ preserves colors, i.e., $\varepsilon(a) = \varepsilon(\phi(a))$ for any $a \in P$, and

\item for all $a <_P b$, $\phi(a) <_{\mathbb{P}_m} \phi(b)$ if and only if $a<_{\mathbb{P}_m} b$.
\end{enumerate}
See Figures \ref{fig:equiv} and \ref{fig:inequiv}. Let $\mathcal{P}^{(m)}_n$ denote the vector space over $\mathbb{Q}$ with basis consisting of all $(m,n)$-posets, modulo equivalence of labelings, and define
\[\mathcal{P}^{(m)} = \bigoplus_{n\geq 0} \mathcal{P}^{(m)}_n.\]
We will define a product $\sqcup_m$ and coproduct $\delta_m$ that make $\mathcal{P}^{(m)}$ into a graded Hopf algebra.

If $P$ is an $(m,n)$-poset and $Q$ is an $(m,r)$-poset then let $P \sqcup_m Q$ be the $(m,n+r)$-poset defined as follows. If as sets $P$ and $Q$ have any elements of the same absolute value, then replace $Q$ by another $(m,r)$-poset that is label-equivalent to $Q$ and whose elements have absolute values distinct from those of $P$. Again, this is easy to do since $P$ and $Q$ are finite sets. Now let $P \sqcup_m Q$ be the poset formed by taking the union of $P$ and $Q$ as posets. We have
\[ |P| \cap |Q| = \emptyset, \text{ and } \]
\[ x <_{P \sqcup_m Q} y \iff x <_P y \text{ or } x <_Q y.\]
Clearly $P \sqcup_m Q \sim Q \sqcup_m P$, and if both $P \sim P'$ and $Q \sim Q'$ then $P \sqcup_m Q \sim P'\sqcup_m Q'$. Thus $\sqcup_m$ can be viewed as a commutative product, giving $\mathcal{P}^{(m)}$ the structure of a graded algebra. The identity element is the empty poset $\emptyset$.

Order ideals of an $(m,n)$-poset are defined just as in the ordinary case. Namely, $I$ is an order ideal of an $(m,n)$-poset $P$ if $x \in I$ and $y <_P x$ implies $y \in I$. Let $\mathcal{I}(P)$ denote the set of order ideals of $P$. The coproduct is defined analogously to the ordinary case: $\delta_m : \mathcal{P}^{(m)} \to \mathcal{P}^{(m)} \otimes \mathcal{P}^{(m)}$ by
\[\delta_m(P) = \sum_{I \in \mathcal{I}(P)} I \otimes (P\setminus I).\] The counit is again projection onto $\mathcal{P}^{(m)}_0 = \mathbb{Q}$.

If $P$ and $Q$ are $m$-colored posets, then
\[\mathcal{I}(P\sqcup_m Q) = \{ I \sqcup_m J ~|~ I \in \mathcal{I}(P) \text{ and } J \in \mathcal{I}(Q)\}.\]
Therefore
\begin{align*} \delta_m(P \sqcup_m Q) & = \sum_{I \in \mathcal{I}(P) \atop J \in \mathcal{I}(Q)} (I \sqcup_m J) \otimes ((P\setminus I) \sqcup_m (Q\setminus J))\\
& = \delta_m(P) \sqcup_m \delta_m(Q).
\end{align*}
So $\delta_m$ is a morphism of algebra and it follows that $\mathcal{P}^{(m)}$ is a graded Hopf algebra with product $\sqcup_m$ and coproduct $\delta_m$. The antipode formula is exactly the same as for ordinary labeled posets, given by \eqref{eq:Pantipode}.

Before continuing, we note that $(m,n)$-chains are naturally related to $m$-colored permutations. Recall that $m$-colored permutations of $n$ are elements of the wreath product $\mathcal{C}_m \wr \mathfrak{S}_n$. We write an element $\pi = \pi_1 \pi_2 \cdots \pi_n \in \mathcal{C}_m \wr \mathfrak{S}_n$ as a word in the letters $[n]_m$ such that $|\pi| = |\pi_1| |\pi_2| \cdots |\pi_n|$ is an ordinary permutation in $\mathfrak{S}_n$. Colored permutations are $(m,n)$-chains under the order inherited from the subscripts: $\pi_1 <_{\pi} \pi_2 <_{\pi} \cdots <_{\pi} \pi_n$. Similarly, the set of linear extensions of an $(m,n)$-poset $P$ is naturally identified with a set of colored permutations.

\begin{remark}
Just as there is a Hopf algebra of permutations, there is a colored permutation Hopf algebra. See \cite{BaumannHohlweg} and \cite{BergeronHohlweg}.
\end{remark}

\subsection{Colored quasisymmetric functions}

The colored quasisymmetric functions are first due to S. Poirier \cite{Poirier}, though N. Bergeron and C. Hohlweg \cite{BergeronHohlweg} have contributed most of the results that we use here.

\begin{figure}
\begin{picture}(90,95)(40,20) \thicklines
   \put(30,50){\line(1,-1){20}}
\put(30,70){\line(1,1){20}}
\put(90,70){\line(-1,1){20}}\put(90,50){\line(-1,-1){20}}
\put(53,20){(3)} \put(15,55){(2,1)} \put(85,55){(1,2)} \put(43,95){(1,1,1)}
\end{picture}
\qquad
\begin{picture}(90,95)(20,20) \thicklines
   \put(30,50){\line(1,-1){20}}
\put(30,70){\line(1,1){20}}
\put(90,70){\line(-1,1){20}}\put(90,50){\line(-1,-1){20}}
\put(50,20){({\color{blue}$\omega 3$})} \put(10,55){({\color{blue}$\omega 2$},{\color{blue}$\omega 1$})} \put(80,55){({\color{blue}$\omega 1$},{\color{blue}$ \omega 2$})} \put(33,95){({\color{blue}$\omega 1$},{\color{blue}$ \omega 1$},{\color{blue}$ \omega 1$})}
\end{picture}
\\
\begin{picture}(60,75)(30,10) \thicklines
   \put(30,30){\line(0,1){20}}
\put(15,20){(2,{\color{blue}$\omega 1$})} \put(12,55){(1,1,{\color{blue}$\omega 1$})}
\end{picture}
\begin{picture}(60,75)(30,10) \thicklines
   \put(30,30){\line(0,1){20}}
\put(15,20){({\color{blue}$\omega 2$},1)} \put(7,55){({\color{blue}$\omega 1$},{\color{blue}$\omega 1$},1)}
\end{picture}
\begin{picture}(60,75)(30,10) \thicklines
   \put(30,30){\line(0,1){20}}
\put(15,20){(1,{\color{blue}$\omega 2$})} \put(9,55){(1,{\color{blue}$\omega 1$},{\color{blue}$\omega 1$})}
\end{picture}
\begin{picture}(60,75)(30,10) \thicklines
   \put(30,30){\line(0,1){20}}
\put(15,20){({\color{blue}$\omega 1$},2)} \put(12,55){({\color{blue}$\omega 1$},1,1)}
\end{picture}
\begin{picture}(60,75)(30,10)
\put(5,35){({\color{blue}$\omega 1$},1,{\color{blue}$\omega 1$}) \quad (1,{\color{blue}$\omega 1$},1)}
\end{picture}
\caption{The partial order on 2-colored compositions of 3.}\label{fig:mcomps}
\end{figure}
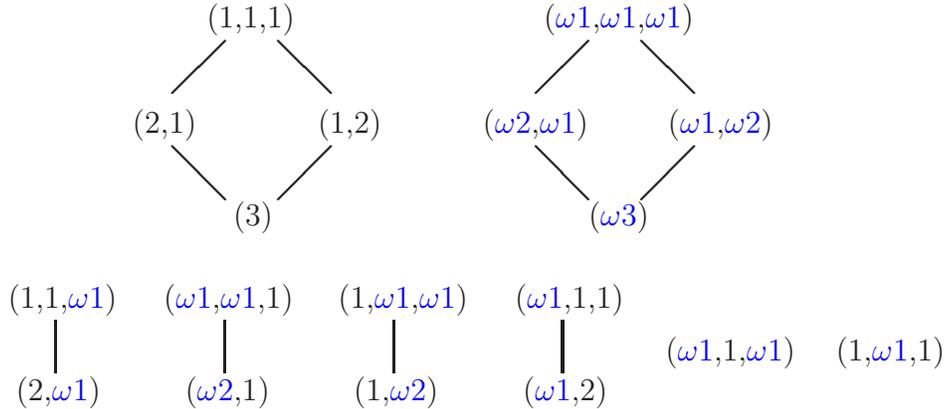

Before defining colored quasisymmetric functions, it is helpful to describe their indexing set: the colored compositions. Fix a positive integer $m$ as the number of colors. An $m$-colored composition is an
ordered tuple of colored positive integers, say $\alpha = (\omega^{j_1} \alpha_1, \omega^{j_2} \alpha_2, \ldots, \omega^{j_k} \alpha_k)$, where $\omega$ is a primitive $m$th root of unity. We write $\alpha \models_m n$ if $|\alpha|:= \alpha_1 + \alpha_2 + \cdots + \alpha_k = n$. For example, if $m =3$, then $\alpha = ( \omega 2, 1, \omega^2 1, 3 )$ is a 3-colored composition of $2+1+1+3 = 7$. The partial order on $m$-colored compositions of $n$ is again by refinement, so that cover relations take the form:\[ (\omega^{j_1} \alpha_1, \ldots, \omega^{j_i} (\alpha_i + \alpha_{i+1}), \ldots, \omega^{j_k} \alpha_k) < (\omega^{j_1} \alpha_1, \ldots, \omega^{j_i} \alpha_i, \omega^{j_i} \alpha_{i+1}, \ldots, \omega^{j_k} \alpha_k). \] Notice in general that the poset of colored compositions of $n$ is disconnected. See Figure \ref{fig:mcomps} for example.

Colored quasisymmetric functions are a natural generalization of quasisymmetric functions to an alphabet with several colors for its letters. For fixed $m$, we consider formal series in the alphabet \[ X^{(m)} := \{ x_{1,0}, x_{1,1}, \ldots, x_{1,m-1}, x_{2,0}, x_{2,1}, \ldots, x_{2,m-1}, \ldots \},\] (so the second subscript corresponds to color) with the same quasisymmetric property. Namely, an $m$-colored quasisymmetric function $Q(X^{(m)})$ is a formal series of bounded degree such that for any $m$-colored composition $\alpha= (\omega^{j_1} \alpha_1, \omega^{j_2} \alpha_2, \ldots, \omega^{j_k} \alpha_k)$, the coefficient of $x_{i_1,j_1}^{\alpha_1} x_{i_2,j_2}^{\alpha_2} \cdots x_{i_k,j_k}^{\alpha_k}$ is the same for all $(i_1,j_1) < (i_2,j_2) < \cdots < (i_k,j_k)$ in the lexicographic order on $\mathbb{P} \times [0,m-1]$. Intuitively, the letters are colored the same as the parts of $\alpha$ and the subscript pairs must be strictly increasing. Let $\Q^{(m)}_n$ denote the set of all $m$-colored quasisymmetric functions homogeneous of degree $n$. Then $\Q^{(m)} := \bigoplus_{n \geq 0} \Q^{(m)}_n$ denotes the graded ring of all $m$-colored quasisymmetric functions, where $\Q^{(m)}_0 = \mathbb{Q}$.

The monomial basis for $\Q^{(m)}_n$ is, for any composition $\alpha = (\omega^{j_1} \alpha_1, \omega^{j_2} \alpha_2, \ldots, \omega^{j_k} \alpha_k) \models_m n$,
\[ M^{(m)}_{\alpha} := \sum_{(i_1,j_1) < (i_2,j_2) < \cdots < (i_k,j_k)} x_{i_1, j_1}^{\alpha_1} x_{i_2, j_2}^{\alpha_2} \cdots x_{i_k, j_k}^{\alpha_k},\] where we remember that the $j_s$, $s=1,\ldots,k$, are fixed by $\alpha$.
There are $m(m+1)^{n-1}$ $m$-colored compositions of $n$, so the graded component $\Q^{(m)}_n$ has dimension $m(m+1)^{n-1}$ as a vector space. The fundamental basis can be written \[ F^{(m)}_{\alpha} := \sum_{ \alpha \leq \beta } M^{(m)}_{\beta}.\] Let $\alpha = (\omega^{j_1} \alpha_1, \ldots, \omega^{j_k} \alpha_k)$. If $s = \alpha_1 + \cdots + \alpha_r + h$, $1 \leq h \leq \alpha_{r+1}$, then define $j'_s = j_{r+1}$, the color of part $\alpha_{r+1}$. We can also write \[F_{\alpha} = \sum_{ (i_1,j'_1) \leq (i_2,j'_2) \leq \cdots \leq (i_n,j'_n) } x_{i_1,j'_1} x_{i_2,j'_2} \cdots x_{i_n,j'_n}, \] where if $s = \alpha_1 + \cdots + \alpha_r$, $r=1,\ldots,l(\alpha)$ then $(i_s,j'_s) < (i_{s+1},j'_{s+1})$ (and the colors are fixed by $\alpha$). For example, \[ F^{(2)}_{1\bar{2}\bar{1}} = \sum_{i \leq j \leq k < l} x_i \bar{x}_j \bar{x}_k \bar{x}_l = M^{(2)}_{1\bar{2}\bar{1}} + M^{(2)}_{1\bar{1}\bar{1}\bar{1}}\] and \[ F^{(3)}_{2\bar{\bar{ 1}} \bar{2}} = \sum_{i \leq j \leq k < l \leq m} x_i x_j \bar{\bar{x}}_k \bar{x}_l \bar{x}_m = M^{(3)}_{2\bar{\bar{1}}\bar{2}} + M^{(3)}_{11\bar{\bar{1}}\bar{2}} + M^{(3)}_{2\bar{\bar{1}}\bar{1}\bar{1}} + M^{(3)}_{11\bar{\bar{1}}\bar{1}\bar{1}},\] where for ease of notation $x_i = x_{i,0}$, $\bar{x}_i = x_{i,1}$, and $\bar{\bar{x}}_i = x_{i,2}$.

Bergeron and Hohlweg \cite{BergeronHohlweg} have shown that $\Q^{(m)}$ is a Hopf algebra.  The product is multiplication of formal power series in $\mathbb{Q}[[X^{(m)}]]$, the counit maps a function to its constant term, and the coproduct, \[\Delta_m : \Q^{(m)} \to \Q^{(m)} \otimes \Q^{(m)},\] is written analogously to the ordinary case for the monomials:
\[ \Delta_m( M^{(m)}_{\alpha} ) = \sum_{ \beta \gamma = \alpha} M^{(m)}_{\beta} \otimes M^{(m)}_{\gamma},\] where $\beta \gamma$ is now the concatenation of colored compositions.

Again the coproduct can be understood in the following way. We let $\mathbb{P}_m + \mathbb{P}_m'$ denote the totally ordered set \[ 1 < \omega 1 < \cdots < \omega^{m-1}1 < 2 < \omega 2 < \cdots < \omega^{m-1}2 < \cdots \qquad\qquad\qquad \] \[ \qquad\qquad\qquad < 1' < \omega 1' < \cdots < \omega^{m-1}1' < 2' < \omega 2' < \cdots < \omega^{m-1}2' < \cdots.\] Then $X^{(m)}+Y^{(m)}$ denotes the set of commuting variables $\{ x_{s,j} : \omega^j s \in \mathbb{P}_m +\mathbb{P}_m'\}$, with the convention $x_{i',j} = y_{i,j}$. For any colored quasisymmetric function $Q$, the coproduct is equivalent to the map $Q(X^{(m)}) \mapsto Q(X^{(m)}+Y^{(m)})$. For example,
\begin{align*}
M^{(2)}_{\bar{2}1}(X^{(2)}+Y^{(2)}) & = \sum_{\omega s < t \in \mathbb{P}_m + \mathbb{P}_m'} x^2_{s,1} x_{t,0} = \sum_{ i < j} x_{i,1}^2 x_{j,0} + \sum_{i,j} x_{i,1}^2 y_{j,0} + \sum_{i < j} y_{i,1}^2 y_{j,0} \\
 & = M^{(2)}_{\bar{2}1}(X^{(2)}) + M^{(2)}_{\bar{2}}(X^{(2)})M^{(2)}_1(Y^{(2)}) + M^{(2)}_{\bar{2}1}(Y^{(2)}),
\end{align*}
and thus \[ \Delta_2(M^{(2)}_{\bar{2}1}) = M^{(2)}_{\bar{2}1} \otimes 1 + M^{(2)}_{\bar{2}} \otimes M^{(2)}_1 + 1 \otimes M^{(2)}_{\bar{2}1}.\]

In Section \ref{sec:colorantipode} we will derive simple antipode formulas for $\Q^{(m)}$.

\subsection{Colored $P$-partitions}

The fundamental colored quasisymmetric functions are naturally related to both descent sets of colored permutations as well as the sequence of colors in a given permutation. These two pieces of information are encoded in colored compositions as we explain shortly. As a bridge between the descent classes of colored permutations and the colored quasisymmetric functions, we will define \emph{colored $P$-partitions}.

Intuitively, the \emph{colored descent composition} of $\pi \in \mathcal{C}_m \wr \mathfrak{S}_n$, $C(\pi)$, is the ordered tuple listing the lengths of increasing runs of $\pi$ with constant color, where we record not only the length of such a run, but also its color. (The notation is unambiguous since the ordinary descent composition is simply a special case.) An example should make the notion clear. If we have two colors (indicated with a bar), let \[\pi = 1\,\bar{2}\bar{3}\,4\,\bar{8}\,\bar{5}\,7\,6.\]
Then \[C(\pi) = (1,\bar{2},1,\bar{1},\bar{1},1,1).\] Again we use $\Des(\pi)$ to denote the set of descents, but now we only count descents between letters of the same color. For example, $\Des( 1\,\bar{2}\bar{3}\,4\,\bar{8}\,\bar{5}\,7\,6)) = \{5,7\}$.  (Warning: Our definition of descent in the case of $2$-colored permutations is different from the standard notion of descent in a Coxeter group of type $B$.)

Now we present colored $P$-partitions to capture this feature. Intuitively, these maps are just color-preserving $P$-partitions. Let $P$ be any labeled $m$-colored poset and let $S$ be any totally ordered countable set.

\begin{defn}
The map $f: P \to S_m$ is a \emph{$m$-colored $P$-partition} if it satisfies the following conditions:
\begin{enumerate}

\item $\varepsilon(i) = \varepsilon(f(i))$ for every $i \in P$

\item $f(i) \leq f(j)$ if $i <_P j$

\item $f(i) < f(j)$ if $i <_P j$ and $i >_{\mathbb{P}_m} j$.

\end{enumerate}
\end{defn}

We let $\mathcal{A}^{(m)}(P)$ (or $\mathcal{A}^{(m)}(P;S_m)$) denote the set of all $m$-colored $P$-partitions, and encode this set in the generating function \[ \Gamma_m(P) := \sum_{ f \in \mathcal{A}^{(m)}(P)} \prod_{i \in P} x_{f(i)},\] where $x_{\omega^j s} = x_{s,j}$. If we take $S$ to be the set of positive integers, then it follows that $\Gamma_m(P) \in \Q^{(m)}_n$. We have that if $P \sim Q$, then $\mathcal{A}^{(m)}(P) = \mathcal{A}^{(m)}(Q)$ and $\Gamma_m(P) = \Gamma_m(Q)$.

For chains,
\begin{align*}
 \mathcal{A}^{(m)}(\pi) = \{ f: \pi \to S_m  \mid & \,\, f(\pi_1) \leq f(\pi_2) \leq \cdots \leq f(\pi_n) \\
 & \mbox{ with } \varepsilon(\pi_i) = \varepsilon(f(\pi_i)) \\
 & \mbox{ and } k \in \Des(\pi) \Rightarrow f(\pi_k) < f(\pi_{k+1})\},
\end{align*}
and if $\varepsilon(\pi_k) = j_k$, the color of $\pi_k$, then
\begin{align*}\label{eq:FmDes} \Gamma_m(\pi) &= \sum_{ \substack{ (i_1,j_1) \leq (i_2,j_2) \leq \cdots \leq (i_n,j_n) \\ k \in \Des(\pi) \Rightarrow (i_k,j_k) < (i_{k+1},j_{k+1}) }} x_{i_1,j_1} x_{i_2,j_2} \cdots x_{i_n,j_n} \\
 &= F^{(m)}_{C(\pi)},
\end{align*}
so the distinct generating functions for the colored $\pi$-partitions of permutations $\pi \in \mathcal{C}_m \wr \mathfrak{S}_n$ form a basis for $\Q^{(m)}_n$.

\begin{remark}
The colored analog of Solomon's descent algebra is the Mantaci-Reutenauer algebra, $MR^{(m)}_n$. The Mantaci-Reutenauer algebra is a subalgebra of the group algebra of $\mathcal{C}_m \wr \mathfrak{S}_n$ given by the span of sums of permutations with common colored descent composition. In \cite{Petersen1} a $P$-partition argument is used to show duality between $MR^{(2)}_n$ and $\Q^{(2)}_n$, but the reasoning should work in general.
\end{remark}

We will now derive some useful facts about colored $P$-partitions by analogy with those of Section \ref{sec:ppartitions}.

\begin{lem}
The set of all $m$-colored $P$-partitions is the disjoint union of the sets of $m$-colored $\pi$-partitions of all linear extensions $\pi$ of $P$: \[\mathcal{A}^{(m)}(P) = \coprod_{\pi \in \mathcal{L}(P)}
\mathcal{A}^{(m)}(\pi).\]
\end{lem}

\begin{proof}
The proof follows by induction on the number of incomparable pairs of elements of $P$. If there are no incomparable pairs in $P$, it is already a chain and the statement is trivially true. Now suppose $i$ and $j$ are incomparable in $P$. Let $P_{ij}$ denote the poset formed from $P$ by introducing the extra relation $i < j$. It is easy to see that $\mathcal{A}^{(m)}(P) = \mathcal{A}^{(m)}(P_{ij}) \coprod \mathcal{A}^{(m)}(P_{ji})$. We continue to split these posets (each with strictly fewer incomparable pairs)
until we have a collection of totally ordered chains corresponding to distinct linear
extensions of $P$.
\end{proof}

\begin{cor}
For any $m$-colored poset $P$,
\begin{equation}\label{eq:mFTPP}
\Gamma_m(P) = \sum_{\pi \in \mathcal{L}(P)} \Gamma_m(\pi).
\end{equation}
\end{cor}

\begin{lem}
The set of $m$-colored $(P \sqcup_m Q)$-partitions is in bijection with the cartesian product of the set of $m$-colored $P$-partitions and the set of $m$-colored $Q$-partitions: \[ \mathcal{A}^{(m)}(P\sqcup_m Q) \leftrightarrow \mathcal{A}^{(m)}(P) \times \mathcal{A}^{(m)}(Q).\]
\end{lem}

\begin{proof}
Given a $P$-partition $f$ and a $Q$-partition $g$ it is obvious that the map $h : P \sqcup_m Q \to S_m$ defined by \[ h(i) =
\begin{cases}
f(i) & \mbox{if } i \in P \\
g(i) & \mbox{if } i \in Q
\end{cases}\] is a $(P \sqcup_m Q)$-partition. Because we can assume the elements of $P$ and $Q$ to be distinctly labeled, this process is reversible.
\end{proof}

\begin{cor}
For any $m$-colored posets $P$ and $Q$,
\begin{equation}\label{eq:mGammamult}
\Gamma_m(P \sqcup_m Q) = \Gamma_m(P)\Gamma_m(Q).
\end{equation}
\end{cor}

Equation \eqref{eq:mGammamult} has the following formula as a special case:
\begin{equation}\label{eq:FMmult}
F^{(m)}_{C(\sigma)} F^{(m)}_{C(\tau)} = \sum_{ \pi \in \Sh(\sigma,\tau)} F^{(m)}_{C(\pi)},
\end{equation}
where the sum is over all shuffles of chains $\sigma$ and $\tau$.

\subsection{Viewing $\Q^{(m)}$ as the homomorphic image of $\mathcal{P}^{(m)}$}

We now use the colored $P$-partition generating function to define the map $\Gamma_m : \mathcal{P}^{(m)} \to \Q^{(m)}$ by $P \mapsto \Gamma_m(P)$. Note that $\Gamma_m(P) = \Gamma_m(Q)$ whenever $P\sim Q$, so $\Gamma_m$ is a well-defined linear map from $\mathcal{P}^{(m)}$ to $\Q^{(m)}$. Formula \eqref{eq:mGammamult} establishes that $\Gamma_m$ is an algebra homomorphism, and below we will demonstrate that $\Gamma_m$ is a coalgebra homomorphism.

Recall that $\mathbb{P}_m + \mathbb{P}_m'$ denotes the set \[ 1 < \omega 1 < \cdots < \omega^{m-1}1 < 2 < \omega 2 < \cdots < \omega^{m-1}2 < \cdots \qquad\qquad\qquad \] \[ \qquad\qquad\qquad < 1' < \omega 1' < \cdots < \omega^{m-1}1' < 2' < \omega 2' < \cdots < \omega^{m-1}2' < \cdots.\]

\begin{lem}\label{lem:colorcoprod}
For any $m$-colored poset $P$ we have \[ \mathcal{A}^{(m)}(P; \mathbb{P}_m + \mathbb{P}_m') \leftrightarrow \coprod_{I \in \mathcal{I}(P)} \mathcal{A}^{(m)}(I; \mathbb{P}_m) \times \mathcal{A}^{(m)}(P\setminus I; \mathbb{P}_m').\]
\end{lem}

The proof of Lemma \ref{lem:colorcoprod} is identical in structure to that of Lemma \ref{lem:coprod}.

\begin{cor}
For any $m$-colored poset $P$ we have
\begin{equation}\label{eq:colorgammacoprod}
\Gamma_m(P)(X^{(m)}+Y^{(m)}) = \sum_{I \in \mathcal{I}(P)} \Gamma_m(I)(X^{(m)})\Gamma(P\setminus I)(Y^{(m)}).
\end{equation}
\end{cor}

Equation \eqref{eq:colorgammacoprod} gives \[\Delta_m(\Gamma_m(P)) = \sum_{I \in \mathcal{I}(P)} \Gamma_m(I) \otimes \Gamma_m(P\setminus I) = \Gamma_m \otimes \Gamma_m( \delta_m(P)), \] which shows that $\Gamma_m$ is a coalgebra homomorphism, and hence is a Hopf homomorphism.

\begin{thm}
The map $\Gamma_m$ is a homomorphism of Hopf algebras. Moreover, \[\Gamma_m(\mathcal{P}^{(m)}) = \Q^{(m)}. \]
\end{thm}

We also get the following nice formula for the coproduct in terms of the fundamental basis by thinking of the order ideals of a chain. For a colored permutation $\pi = \pi_1 \cdots \pi_n$,
\begin{equation}\label{eq:FMcoproduct}
 \Delta_m(F^{(m)}_{C(\pi)}) = \sum_{i = 0}^n F^{(m)}_{C(\pi_1 \cdots \pi_i)} \otimes F^{(m)}_{C(\pi_{i+1} \cdots \pi_n)}.
\end{equation}

\subsection{Antipode formulas and cycloribbon diagrams}\label{sec:colorantipode}

In this section we give formulas for the antipode $S_m: \Q^{(m)} \to \Q^{(m)}$ on the $M$ and $F$ bases. The formulas are straightforward generalizations of \eqref{eq:antipodeM} and \eqref{eq:Fantipode}.

\begin{prp}\label{prp:colorantipodeM}
For any $m$-colored composition $\alpha$, the antipode $S_m$ satisfies
\[S_m(M^{(m)}_{\alpha}) = (-1)^{l(\alpha)}\sum_{\beta \leq \alpha} M^{(m)}_{\overleftarrow{\beta}},\]
where $\overleftarrow{\beta} = (\beta_k, \beta_{k-1}, \ldots, \beta_1)$, the $m$-colored composition formed by writing the parts of $\beta$ backwards.
\end{prp}

 \begin{prp}\label{prp:FMantipode}
For any $m$-colored permutation $\pi = \pi_1 \cdots \pi_n$,
\begin{equation}\label{eq:FMantipode}
S_m(F^{(m)}_{C(\pi)}) = (-1)^n F^{(m)}_{C(\pi_n \pi_{n-1} \cdots \pi_1)}
\end{equation}
 \end{prp}

\begin{proof}[Proof of Proposition \ref{prp:colorantipodeM}]
This proposition can be proved by working in the dual Hopf algebra of $\Q^{(m)}$, as in \cite[Corollary~2.3]{MalvenutoReutenauer}, but we will present an induction proof following the same lines of reasoning used in \cite[Proposition~3.4]{Ehrenborg}.

If $l(\alpha) = 0$ or 1, the formula is obviously correct. Let $l(\alpha) = k$ and suppose the formula holds for all shorter compositions. Using the inductive formula for the antipode \eqref{eq:antipode}, we have
\begin{align}
S_m(M^{(m)}_{\alpha}) &= -M^{(m)}_{ \alpha} - \sum_{i=1}^{k-1} S_m(M^{(m)}_{(\alpha_1,\ldots,\alpha_i)})M^{(m)}_{(\alpha_{i+1},\ldots,\alpha_k)} \nonumber \\
&= -M^{(m)}_{ \alpha} - \sum_{i=1}^{k-1} (-1)^i \sum_{\beta \leq (\alpha_i,\ldots,\alpha_1) } M^{(m)}_{\beta} M^{(m)}_{(\alpha_{i+1},\ldots,\alpha_k)}.\label{eq:expandS}
\end{align}

Any $\beta \leq (\alpha_i, \ldots, \alpha_1)$ will have a first entry of the form $\alpha_i + \alpha_{i-1} + \cdots + \alpha_j$ where $i \geq j$ and $\alpha_i, \alpha_{i-1}, \ldots, \alpha_j$ all have the same color $\varepsilon$. When expanding the product $M^{(m)}_{\beta} M^{(m)}_{(\alpha_{i+1},\ldots, \alpha_k)}$ into monomials, the first part of each composition that appears will be one of the following types:
\begin{enumerate}
\item $\alpha_i + \alpha_{i-1} + \cdots + \alpha_j$ as above,

\item $\alpha_i + \alpha_{i-1} + \cdots + \alpha_j + \alpha_{i+1}$ (only possible if $\varepsilon$ is also the color of $\alpha_{i+1}$), or

\item $\alpha_{i+1}$.
\end{enumerate}
We say that a monomial appearing in the sum \eqref{eq:expandS} has type $i$ if the first part of its composition is as in the first case. If its composition is as in the second or third cases, we say that the monomial has type $i+1$.

Notice that any monomial of type $i$, $i=1,\ldots,k-1$ that appears in this expression for $S_m(M^{(m)}_{\alpha})$ will occur in both the $i$th and $(i-1)$st term of \eqref{eq:expandS}, and with opposite sign (understanding that $-M^{(m)}_{\alpha}$ is the 0th term). Hence the only terms that don't cancel in \eqref{eq:expandS} are the monomials of type $k$, which appear in the $(k-1)$st term with sign $(-1)^k$.

Now observe that each monomial $M^{(m)}_{\gamma}$ of type $k$ has $\overleftarrow{\alpha}$ as a refinement, i.e., $\gamma \leq \overleftarrow{\alpha}$. On the other hand, if $\overleftarrow{\alpha}$ refines $\gamma$, then $\gamma$ must be of type $k$.
\end{proof}

The following proof is also by induction.

\begin{proof}[Proof of Proposition \ref{prp:FMantipode}]
If $n = 0$ or 1 the proposition holds. Assume that $n\ge 1$ and suppose by induction the proposition holds for $m$-colored permutations with fewer than $n$ letters. For $1\le i\le n-1$, let $\pi_{\le i} = \pi_i \pi_{i-1} \cdots \pi_1$  and $\pi_{> i} = \pi_{i+1} \cdots \pi_n$. Using the recursive antipode formula \eqref{eq:antipode}, together with \eqref{eq:FMmult} and \eqref{eq:FMcoproduct}, we get
\[ S_m(F^{(m)}_{C(\pi)}) = -F^{(m)}_{C(\pi)} - \sum_{i = 1}^{n-1} (-1)^{i+1} \sum_{\sigma \in \Sh(\pi_{\le i}, \pi_{>i})}  F^{(m)}_{C(\sigma)}.\]
Given two disjoint $m$-colored permutations $\sigma=\sigma_1 \cdots \sigma_k$ and $\tau=\tau_1 \cdots \tau_l$, let $\Sh_1(\sigma,\tau)$ denote the set of shuffles of $\sigma$ and $\tau$ that begin with $\sigma_1$ and let $\Sh_2(\sigma,\tau)$ denote the set of shuffles of $\sigma$ and $\tau$ that begin with $\tau_1$. Thus $\Sh(\sigma,\tau)$ is the disjoint union of $\Sh_1(\sigma,\tau)$ and $\Sh_2(\sigma,\tau)$, allowing us to write
\begin{align*} S_m(F^{(m)}_{C(\pi)}) & =  \sum_{i = 0}^{n-2} (-1)^{i+1} \left(\sum_{\sigma \in \Sh_2(\pi_{\le i}, \pi_{> i})}  F^{(m)}_{C(\sigma)}  - \sum_{\sigma \in \Sh_1(\pi_{\le i+1}, \pi_{> i+1})}  F^{(m)}_{C(\sigma)} \right) \\
& \quad + (-1)^n \sum_{\sigma\in\Sh_2(\pi_{\le n-1}, \pi_{> n-1})} F^{(m)}_{C(\sigma)},
\end{align*}
where we set $\pi_{\le 0}  = \emptyset$. Since clearly $\Sh_2(\pi_{\le i}, \pi_{> i}) = \Sh_1(\pi_{\le i+1}, \pi_{> i+1})$ for $0 \le i \le n-2$, the first sum vanishes, leaving us with the second sum, which is none other than \[(-1)^n F^{(m)}_{C(\pi_n \pi_{n-1} \cdots \pi_1)}\] as claimed.
\end{proof}

There is a straightforward way to determine $S(F^{(m)}_\alpha)$ without referring to any underlying permutation. This involves encoding colored compositions by cycloribbon diagrams, as introduced by Hivert, Novelli, and Thibon \cite{HivertNovelliThibon}. First recall that a ribbon diagram consists of rows of empty squares, with the last square of each row directly above the first square of the next row. A \emph{cycloribbon diagram} is a ribbon diagram in which each square is filled in with a color such that the colors in each row are weakly increasing and the colors in each column are weakly decreasing.

\begin{figure}
 \begin{picture}(140,120){
\put(0,48){\framebox(12,12){$0$}} \put(12,48){\framebox(12,12){$2$}}
\put(12,36){\framebox(12,12){$1$}}\put(24,36){\framebox(12,12){$1$}}
\put(24,24){\framebox(12,12){$1$}}
\put(36,24){\framebox(12,12){$1$}}
\put(48,24){\framebox(12,12){$1$}}\put(60,24){\framebox(12,12){$2$}}
\put(60,12){\framebox(12,12){$2$}}\put(72,12){\framebox(12,12){$2$}}
\put(72,0){\framebox(12,12){$0$}}\put(84,0){\framebox(12,12){$0$}}
\put(96,0){\framebox(12,12){$0$}}\put(108,0){\framebox(12,12){$0$}}
}
\end{picture}
\begin{picture}(60,120){
\put(0,108){\framebox(12,12){$0$}} \put(0,96){\framebox(12,12){$0$}}
\put(0,84){\framebox(12,12){$0$}}\put(0,72){\framebox(12,12){$0$}}
\put(12,72){\framebox(12,12){$2$}}
\put(12,60){\framebox(12,12){$2$}}
\put(24,60){\framebox(12,12){$2$}}\put(24,48){\framebox(12,12){$1$}}
\put(24,36){\framebox(12,12){$1$}}\put(24,24){\framebox(12,12){$1$}}
\put(36,24){\framebox(12,12){$1$}}\put(36,12){\framebox(12,12){$1$}}
\put(48,12){\framebox(12,12){$2$}}\put(48,0){\framebox(12,12){$0$}}
}
\end{picture}
\caption{Cycloribbons corresponding to $\alpha=(1, \bar{\bar{1}}, \bar{2}, \bar{3}, \bar{\bar{1}}, \bar{\bar{2}}, 4)$ and $\widetilde{\alpha} = (1,1,1,1,\bar{\bar{1}},\bar{\bar{2}},\bar{1},\bar{1},\bar{2},\bar{1},\bar{\bar{1}},1)$. }\label{fig:cycloribbons}
\end{figure}
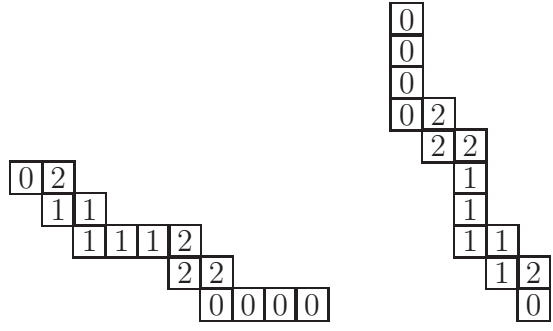

Now given a colored composition $\alpha = (\omega^{j_1}\alpha_1, \ldots, \omega^{j_k}\alpha_k)$, make a row of $\alpha_i$ squares colored $j_i$ for every $i\le k$ and arrange them as follows: if $j_i < j_{i+1}$ then extend row $i$ by adjoining row $i+1$ to the end of row $i$, and if $j_i\ge j_{i+1}$ then put the first square in row $i+1$ directly below the last square in row $i$.
The resulting diagram is a cycloribbon diagram. Conversely, every cycloribbon diagram encodes a unique colored composition.  If the cycloribbon diagram of $\alpha$ is reflected across the line $y=x$ so that rows and columns are switched, then the resulting diagram is again a cycloribbon diagram, and the corresponding colored composition is called the \emph{conjugate} of $\alpha$, denoted $\widetilde{\alpha}$. Figure \ref{fig:cycloribbons} shows cycloribbon diagrams for the $3$-colored composition $\alpha=(1, \bar{\bar{1}}, \bar{2}, \bar{3}, \bar{\bar{1}}, \bar{\bar{2}}, 4)$ (left) and its conjugate $\widetilde{\alpha} = (1,1,1,1,\bar{\bar{1}},\bar{\bar{2}},\bar{1},\bar{1},\bar{2},\bar{1},\bar{\bar{1}},1)$ (right), where the number of bars indicates the color.

It is easy to verify that if $\alpha = C(\pi_1 \cdots \pi_n)$, then $\widetilde{\alpha} = C(\pi_n \pi_{n-1} \cdots \pi_1)$, so by Proposition~\ref{prp:FMantipode} we get the following result:
\begin{cor}
For any $m$-colored composition $\alpha$, the antipode $S_m$ satisfies
\[ S_m(F^{(m)}_\alpha) = (-1)^n F^{(m)}_{\widetilde{\alpha}}.\]
\end{cor}

\subsection{Colored peak functions}

Bergeron and Hohlweg have defined a Hopf subalgebra of $\Q^{(m)}$ analogous to Stembridge's peak algebra. Their first definition is given in terms of the fundamental basis, but we will write it in terms of the monomials by analogy with section \ref{sec:peakA}. For any $\alpha \models_m n$, we can write $\alpha = \alpha_{(1)} \alpha_{(2)} \cdots \alpha_{(k)}$, the concatenation of compositions $\alpha_{(i)}$, such that all the parts of each $\alpha_{(i)}$ have the same color, $\epsilon_i$, and no two consecutive $\alpha_{(i)}$ have the same color. Bergeron and Hohlweg refer to this way of writing a colored composition as the \emph{rainbow decomposition} of $\alpha$. For example, with $m=3$, \[ (\bar{2},\bar{1},1,3,2,\bar{\bar{1}},\bar{2},\bar{4}) = (\bar{2},\bar{1})(1,3,2)(\bar{\bar{1}})(\bar{2},\bar{4}).\] Now if $\alpha = \alpha_{(1)} \alpha_{(2)} \cdots \alpha_{(k)}$, we define
\[K^{(m)}_{\alpha} := \sum_{ \alpha_{(i)} \leq \epsilon_i \beta_{(i)}^* } 2^{\sum l(\beta_{(i)})} M^{(m)}_{\epsilon_1\beta_{(1)} \epsilon_2 \beta_{(2)} \cdots \epsilon_k \beta_{(k)}},\]
where $\epsilon_i \beta_{(i)}$ means that we give each of the parts of $\beta_{(i)}$ the color $\epsilon_i$, and as before $\beta^*$ is the refinement of $\beta$ obtained by replacing, for $i > 1$, every part $\beta_i \geq 2$ with the two parts $(1,\beta_i -1)$. Essentially, this definition is the same as \eqref{eq:KA-M} but we handle the differently colored parts of $\alpha$ independently. For example, \[K^{(2)}_{21\bar{1}} = 2^3 M^{(2)}_{21\bar{1}} + 2^3 M^{(2)}_{12\bar{1}} + 2^4 M^{(2)}_{111\bar{1}}.\]

Consider all $\alpha = \alpha_{(1)} \alpha_{(2)} \cdots, \alpha_{(k)} \models_m n$ such that each $\alpha_{(i)}$ is a peak composition in the color $\epsilon_i$. We call such compositions \emph{$m$-colored peak compositions}. Let $\mathbf{\Pi}^{(m)}_n$ denote the span of the functions $K^{(m)}_{\alpha}$, over all $m$-colored peak compositions $\alpha$. These $K^{(m)}_{\alpha}$ are the $m$-colored peak functions. By their definition in \cite{BergeronHohlweg}, the colored peak functions are $F^{(m)}$-positive and linearly independent, forming a basis for $\mathbf{\Pi}_n$. We define $\mathbf{\Pi}^{(m)} := \bigoplus_{n \geq 0} \mathbf{\Pi}^{(m)}_n$ to be the graded ring of colored peak functions. The dimensions of the graded components are given by ($m$ times the) $m$-Fibonacci numbers, i.e., the Fibonacci polynomials evaluated at $x=m$.

\begin{prp}\label{prp:peakdim}
The number of $m$-colored peak compositions of $n$, $f_{m,n}$, is given by the recurrence \[ f_{m,n} = m f_{m,n-1} + f_{m,n-2},\] with initial values $f_{m,1} = m$ and $f_{m,2} = m^2$.
\end{prp}

\begin{cor}
\[\dim(\mathbf{\Pi}^{(m)}_n) = f_{m,n}.\]
\end{cor}
%

\begin{proof}[Proof of Proposition \ref{prp:peakdim}]

Clearly, $f_{m,1} = m$ since $(\omega^j 1)$ is a peak composition for any $j = 0,1,\ldots,m-1$. It is also easy to check that $f_{m,2} = m^2$ since we have $m$ peak compositions of the form $(\omega^j 2)$ and $m(m-1)$ peak compositions of the form $(\omega^i 1, \omega^j 1)$, $i \neq j$.

For the induction, consider any peak composition $\alpha = (\epsilon_1 \alpha_1, \ldots, \epsilon_k \alpha_k) \models_m n-1$. We form a peak composition of $n$ by adding a part of size 1 to the end of $\alpha$ with any color $\epsilon' \neq \epsilon_k$, or by adding 1 to $\alpha_k$ with the same color: \[(\epsilon_1 \alpha_1, \ldots, \epsilon_k \alpha_k) \mapsto \begin{cases} (\epsilon_1 \alpha_1, \ldots, \epsilon_k (\alpha_k + 1)) \mbox{ or}\\ (\epsilon_1 \alpha_1, \ldots, \epsilon_k \alpha_k, \epsilon' 1) \mbox{ for any color } \epsilon' \neq \epsilon_k.\end{cases} \] The only peak compositions of $n$ that aren't covered by this case are those whose last part is size 1 and the same color as the next-to-last part. We get these from any of the peak compositions of $n-2$ by adding 1 to the last part and adding a new part of size 1 to the end. For $\alpha = (\epsilon_1 \alpha_1, \ldots, \epsilon_k \alpha_k) \models_m n-2$, \[ (\epsilon_1 \alpha_1, \ldots, \epsilon_k \alpha_k) \mapsto (\epsilon_1 \alpha_1, \ldots, \epsilon_k (\alpha_k+1), \epsilon_k 1).\]

\end{proof}

For any colored composition $\alpha$ with rainbow decomposition $\alpha_{(1)} \alpha_{(2)} \cdots, \alpha_{(k)}$, let $\widehat{\alpha}$ denote the composition formed by replacing consecutive 1s of the same color with their sum added to the next part to the right. In other words, apply the ``\,\,$\widehat{\,\,}$\,\," operation to the $\alpha_{(i)}$ and then concatenate: \[ \widehat{\alpha} = \widehat{\alpha_{(1)}} \widehat{\alpha_{(2)}} \cdots \widehat{\alpha_{(k)}}.\] For example, if $\alpha = (3,1,\bar{1},\bar{3},2,\bar{1},\bar{1},1)$, then $\widehat{\alpha} = (3,1,\bar{4},2,\bar{2},1)$. The map $\alpha \mapsto \widehat{\alpha}$ is a surjection from all colored compositions onto colored peak compositions. Bergeron and Hohlweg \cite[Theorem~5.3]{BergeronHohlweg} show that the induced map $\Theta_m : \Q^{(m)} \to \Q^{(m)}$ given by \[ \Theta_m(F^{(m)}_{\alpha}) = K^{(m)}_{\widehat{\alpha}}\] is a Hopf homomorphism whose image is $\mathbf{\Pi}^{(m)}$. We will motivate this result with the theory of colored enriched $P$-partitions.

\subsection{Colored enriched $P$-partitions}

Bergeron and Hohlweg define the colored peak functions as the graded dual of a colored peak subalgebra of the Mantaci-Reutenauer algebra. By definition then, the $K^{(m)}_{\alpha}$, where $\alpha$ is a colored peak composition, encode peak classes of colored permutations. However, we would like to motivate this connection as we have in the earlier cases. Thus in this section we introduce the colored enriched $P$-partitions.

Just as we have a rainbow decomposition for compositions, we also have a rainbow decomposition for colored permutations. We write $\pi = \pi_{(1)} \pi_{(2)} \cdots \pi_{(k)}$, where each $\pi_{(i)}$ is a permutation in which each letter has the same color $\varepsilon_i$ and no two consecutive words have the same color: $\varepsilon_i \neq \varepsilon_{i+1}$. The peak set of a colored permutation we define to be the set of peaks among runs of the same color in $\pi$. We let $\widehat{C}(\pi)$ denote the peak composition of a colored permutation $\pi$, defined as the concatenation of the peak compositions of the $\pi_{(i)}$: \[\widehat{C}(\pi) = \varepsilon_1 \widehat{C}(\pi_{(1)}) \varepsilon_2 \widehat{C}(\pi_{(2)}) \cdots \varepsilon_k \widehat{C}(\pi_{(k)}).\] Notice that the peak composition encodes more information than simply the peak set. For example, if $\pi = \bar{3}\bar{7}\,\bar{2}\bar{5}\,41\,\bar{8}\bar{9}\,\bar{6}$, then $\widehat{C}(\pi) = (\bar{2},\bar{2},2,\bar{2},\bar{1})$ whereas $\Pe(\pi) = \{2,8\}$ says nothing about color.

Colored enriched $P$-partitions can be thought of as color-preserving enriched $P$-partitions. First we need the proper definition for the image set. For any $m$-colored set $S_m$ with total order \[s_1 < \omega s_1 < \cdots < \omega^{m-1} s_1 < s_2 < \omega s_2 < \cdots < \omega^{m-1} s_2 < \cdots,\] define the set $S_m^{\pm}$ to be \[-s_1 < s_1 < -\omega s_1 < \omega s_1 < \cdots < -\omega^{m-1} s_1 < \omega^{m-1} s_1 \qquad\qquad\qquad\] \[ \qquad\qquad\qquad < -s_2 < s_2 < -\omega s_2 < \omega s_2 < \cdots < -\omega^{m-1} s_2 < \omega^{m-1} s_2 < \cdots.\] We write $s \leq^+ t$ in $S_m^{\pm}$ if $s <_{S_m^{\pm}} t$ or if $s = t \in S_m$. Similarly, $s \leq^{-} t$ means $s <_{S_m^{\pm}} t$ or $-s = -t \in S_m$. We should be clear that in this context the minus signs and colors are not meant to interact in any way. Here the minus sign is just another distinguishing feature of an element. From this point forward, the absolute value $|\cdot| : S^{\pm}_m \to S_m$ simply forgets minus signs: $|- \omega^j s_i | = |\omega^j s_i| = \omega^j s_i$.

\begin{defn}
An \emph{$m$-colored enriched $P$-partition} is a map $f: P \to S_m^{\pm}$ such that
\begin{enumerate}
\item $\varepsilon( i ) = \varepsilon( f(i) )$ for every $i \in P$

\item $f(i) \leq^+ f(j)$ if $i <_P j$ and $i <_{\mathbb{P}_m} j$

\item $f(i) \leq^- f(j)$ if $i <_P j$ and $i >_{\mathbb{P}_m} j$.
\end{enumerate}
\end{defn}

We let $\mathcal{E}^{(m)}(P)$ (or $\mathcal{E}^{(m)}(P;S_m^{\pm})$) denote the set of all $m$-colored enriched $P$-partitions. The generating function for $m$-colored enriched $P$-partitions is \[ \Lambda_m(P) := \sum_{ f \in \mathcal{E}^{(m)}(P)} \prod_{i\in P} x_{|f(i)|},\] remembering our use of absolute value here. If we take $S_m = \mathbb{P}_m$ then $\Lambda_m(P) \in \Q^{(m)}_n$. Also note that if $P \sim Q$, then $\mathcal{E}^{(m)}(P) = \mathcal{E}^{(m)}(Q)$ and $\Lambda_m(P) = \Lambda_m(Q)$.

Following \cite[Proposition~2.2]{Stembridge}, we have the following.

\begin{prp}
For any $m$-colored permutation $\pi$, we have \[ \Lambda_m(\pi) = K^{(m)}_{\widehat{C}(\pi)}.\]
\end{prp}

\begin{proof}
Let $\widehat{C}(\pi)$ have rainbow decomposition $\alpha_{(1)}\cdots \alpha_{(k)}$, where sub-composition $\alpha_{(i)}$ has color $\epsilon_i$. We want to show that
\begin{equation}\label{eq:LambdainM}
\Lambda_m(\pi) = \sum_{ \alpha_{(i)} \leq \epsilon_i \beta_{(i)}^* } 2^{\sum l(\beta_{(i)})} M^{(m)}_{\epsilon_1\beta_{(1)} \epsilon_2 \beta_{(2)} \cdots \epsilon_k \beta_{(k)}}.
\end{equation}

Clearly, we can expand $\Lambda_m(\pi)$ as a sum of monomials $M^{(m)}_{\beta}$ with nonnegative coefficients. By the color-preserving property, we need only consider those compositions $\beta$ with the same type of rainbow decomposition as $\widehat{C}(\pi)$. Fix $\beta = (\omega^{j_1} \beta_1,\ldots, \omega^{j_l} \beta_l) = \beta_{(1)} \cdots \beta_{(k)}$. The coefficient of $M^{(m)}_{\beta}$ is equal to the coefficient of $x_{1,j_1}^{\beta_1} x_{2,j_2}^{\beta_2} \cdots x_{l,j_l}^{\beta_{l}}$. This coefficient is equal to the number of colored enriched $\pi$-partitions $f$ such that
\begin{align}
(|f(\pi_1)|, |f(\pi_2)|, & \ldots, |f(\pi_n)|) \nonumber \\
 & = ( \underbrace{\omega^{j_1} 1,\omega^{j_1} 1,\ldots, \omega^{j_1} 1}_{\beta_1}, \underbrace{ \omega^{j_2} 2,\omega^{j_2} 2,\ldots , \omega^{j_2} 2}_{\beta_2}, \ldots, \underbrace{\omega^{j_l} l,\omega^{j_l} l,\ldots, \omega^{j_l} l}_{\beta_{l}} ). \label{eq:tuple}
\end{align}
Notice that within each block of numbers, the colored permutation $\pi$ must have constant color and satisfy one of three conditions. It must either be always increasing, always decreasing, or decreasing then increasing. It cannot be increasing then decreasing, since then it would have a peak. Say there is a peak in position $r$. Then $f(\pi_{r-1}) \leq^{+} f(\pi_r) \leq^{-} f(\pi_{r+1})$, so $|f(\pi_{r-1})| < |f(\pi_{r+1})|$. Therefore the only possible positions for peaks are on the boundary between two consecutive blocks of the same color.

Let $\Pe_{(i)}(\pi) \subset \Pe(\pi)$ be the subset of the peak set corresponding to $\alpha_{(i)}$ and if $B \subset [n-1]$ is the set corresponding to $\beta$, let $B_{(i)} \subset B$ correspond to $\beta_{(i)}$. Saying that peaks occur only on the boundaries is equivalent to saying $\Pe_{(i)} \subset B_{(i)} \cup (B_{(i)} +1)$ (where $S+1 := \{ s+1 : s \in S\}$), written $\alpha_{(i)} \leq \beta_{(i)}^*$ in terms of compositions. This must be true for every $i$, so the coefficient of $M^{(m)}_{\beta}$ is nonzero only if $\alpha_{(i)} \leq \beta_{(i)}^*$ for all $i$. These are precisely the $M^{(m)}_{\beta}$ occurring on the right-hand side of \eqref{eq:LambdainM}.

Now we will determine the coefficient of such an $M^{(m)}_{\beta}$. Within each block of numbers in \eqref{eq:tuple}, we have: \[ ( |f(\pi_r)|, |f(\pi_{r+1})|, \ldots, |f(\pi_{r+\beta_i})| ) = ( \omega^{j_i} i,\omega^{j_i}i, \ldots, \omega^{j_i} i),\] where $r = \beta_1 + \cdots + \beta_{i-1}+1$. We claim there are exactly two possibilities for $f$ in every such block. If $\pi$ is increasing over this interval, then \[f(\pi_r) \leq^{+} f(\pi_{r+1}) \leq^{+} \cdots \leq^{+} f(\pi_{r+\beta_i}),\] so then $f(\pi_r) = \pm \omega^{j_i} i$, and all others equal $+ \omega^{j_i} i$. If $\pi$ is decreasing over the entire interval, then \[f(\pi_r) \leq^{-} f(\pi_{r+1}) \leq^{-} \cdots \leq^{-} f(\pi_{r+\beta_i}),\] so then $f(\pi_{\beta_i}) = \pm \omega^{j_i} i$, and all others equal $-\omega^{j_i}$. The third case has $\pi$ decreasing, then increasing. Suppose $\pi_{s-1} > \pi_s < \pi_{s+1}$ with $r < s < r+\beta_i$. Then \[f(\pi_r) \leq^{-} \cdots \leq^{-} f(\pi_s) \leq^{+} \cdots \leq^{+} f(\pi_{r+\beta_i}),\] so then $f(\pi_s) = \pm \omega^{j_i} i$, everything in the block to its left is equal to $-\omega^{j_i} i$, while everything to its right is $+\omega^{j_i} i$.

In total, there are $l = \sum l(\beta_{(i)}) $ blocks where some choice can be made, and so there are $2^l$ such $f$. This completes the proof.
\end{proof}

\begin{remark}
Bergeron and Hohlweg \cite[Theorem~2.12]{BergeronHohlweg} show the sums of permutations of $m$-colored permutations of $[n]$ with the same colored peak composition form a subalgebra of the Mantaci-Reutenauer algebra: $\mathfrak{P}^{(m)}_n \subset MR^{(m)}_n$. It seems reasonable to expect that the colored enriched $P$-partitions can be employed to show duality between $\mathfrak{P}^{(m)}_n$ and $\mathbf{\Pi}^{(m)}_n$, but this exercise has not been carried out.
\end{remark}

Now we present some useful facts about colored enriched $P$-partitions. The proofs of these statements follow from the same arguments used earlier for similar statements.

\begin{lem}
The set of all $m$-colored enriched $P$-partitions is the disjoint union of the sets of $m$-colored enriched $\pi$-partitions of all linear extensions $\pi$ of $P$: \[\mathcal{E}^{(m)}(P) = \coprod_{\pi \in \mathcal{L}(P)}
\mathcal{E}^{(m)}(\pi).\]
\end{lem}

\begin{cor}
For any $m$-colored poset $P$,
\begin{equation}\label{eq:mLambdasum}
\Lambda_m(P) = \sum_{\pi \in \mathcal{L}(P)} \Lambda_m(\pi).
\end{equation}
\end{cor}

\begin{lem}
The set of $m$-colored enriched $(P \sqcup_m Q)$-partitions is in bijection with the cartesian product of the set of $m$-colored enriched $P$-partitions and the set of $m$-colored enriched $Q$-partitions: \[ \mathcal{E}^{(m)}(P\sqcup_m Q) \leftrightarrow \mathcal{E}^{(m)}(P) \times \mathcal{E}^{(m)}(Q).\]
\end{lem}

\begin{cor}
For any $m$-colored posets $P$ and $Q$,
\begin{equation}\label{eq:mLambdamult}
\Lambda_m(P \sqcup_m Q) = \Lambda_m(P)\Lambda_m(Q).
\end{equation}
\end{cor}

Taking $P$ and $Q$ to be chains in \eqref{eq:mLambdamult} gives
\begin{equation}\label{eq:colorKmult}
K^{(m)}_{\widehat{C}(\sigma)} K^{(m)}_{\widehat{C}(\tau)} = \sum_{ \pi \in \Sh(\sigma,\tau) } K^{(m)}_{\widehat{C}(\pi)},
\end{equation}
the sum over all shuffles of the colored permutations $\sigma$ and $\tau$, i.e., the image of multiplication of chains in $\mathcal{P}^{(m)}$.

\subsection{Viewing $\mathbf{\Pi}^{(m)}$ as the homomorphic image of $\mathcal{P}^{(m)}$}

Define the map $\Lambda_m: \mathcal{P}^{(m)} \to \Q^{(m)}$ by taking a colored poset to its colored enriched generating function. We will show that this map is a Hopf homomorphism and therefore its image, $\mathbf{\Pi}^{(m)}$, is a Hopf algebra. First we note that equation \eqref{eq:mLambdamult} shows that $\Lambda_m$ is an algebra homomorphism.

To see that it is a coalgebra homomorphism, we use the same trick as before. Let $\mathbb{P}_m^{\pm} + (\mathbb{P}_m^{\pm})'$ denote the set \[-1 < 1 < -\omega 1 < \omega 1 < \cdots < -\omega^{m-1} 1 < \omega^{m-1} 1 \qquad\qquad\qquad\qquad\] \[ \qquad\qquad < -2 < 2 < -\omega 2 < \omega 2 < \cdots < -\omega^{m-1} 2 < \omega^{m-1} 2 < \cdots\] \[< -1' < 1' < -\omega 1' < \omega 1' < \cdots < -\omega^{m-1} 1' < \omega^{m-1} 1' \qquad\qquad\] \[ \qquad\qquad\qquad\qquad < -2' < 2' < -\omega 2' < \omega 2' < \cdots < -\omega^{m-1} 2' < \omega^{m-1} 2' < \cdots.\]

The proof of the following is as in the proof of Lemma \ref{lem:coprod}.

\begin{lem}
For an $m$-colored poset $P$ we have \[ \mathcal{E}^{(m)}(P; \mathbb{P}_m^{\pm} + (\mathbb{P}_m^{\pm})') \leftrightarrow \coprod_{I \in \mathcal{I}(P)} \mathcal{E}^{(m)}(I; \mathbb{P}_m^{\pm}) \times \mathcal{E}^{(m)}(P\setminus I; (\mathbb{P}_m^{\pm})').\]
\end{lem}

\begin{cor}\label{cor:mLambdacoprod}
For an $m$-colored poset $P$ we have
\[
\Lambda_m(P)(X^{(m)}+Y^{(m)}) = \sum_{I \in \mathcal{I}(P)} \Lambda_m(I)(X^{(m)})\Lambda_m(P\setminus I)(Y^{(m)}).
\]
\end{cor}

So we can conclude that $\Lambda_m$ is a coalgebra homomorphism, and hence $\Lambda_m$ is a Hopf homomorphism.

\begin{thm}
The map $\Lambda_m$ is a homomorphism of Hopf algebras. Moreover, \[\Lambda_m(\mathcal{P}^{(m)}) = \mathbf{\Pi}^{(m)}. \]
\end{thm}

 As a special case of Corollary \ref{cor:mLambdacoprod} we can give an explicit formula for the coproduct of a peak function in terms of the ideals of a chain $\pi = \pi_1 \cdots \pi_n$:
\begin{equation}\label{eq:mKcoprod}
\Delta_m(K^{(m)}_{\widehat{C}(\pi)}) = \sum_{i = 0}^n K^{(m)}_{\widehat{C}(\pi_1 \cdots \pi_i)} \otimes K^{(m)}_{\widehat{C}(\pi_{i+1} \cdots \pi_n)}.
\end{equation}

By taking formulas \eqref{eq:FMmult} and \eqref{eq:colorKmult} we can see that $\Theta$ is an algebra homomorphism and formulas \eqref{eq:FMcoproduct} and \eqref{eq:mKcoprod} show that it is also a coalgebra homomorphism. Now it is clear that we can write $\Lambda_m$ as $\Lambda_m = \Theta_m \circ \Gamma_m$, the composition of two Hopf homomorphisms. We get the antipode formula
\[ S_m(K^{(m)}_{\widehat{C}(\pi)}) = (-1)^nK^{(m)}_{\widehat{C}((\pi_n, \pi_{n-1}, \ldots, \pi_1))}\]
by \eqref{eq:FMantipode} and the fact that $S_m\circ \Theta_m = \Theta_m \circ S_m$.

\subsection{Theory of $m$-colored combinatorial Hopf algebras}\label{sec:colorHopf}

Bergeron and Hohlweg \cite{BergeronHohlweg} have outlined a theory of colored combinatorial Hopf algebras. We will present our interpretation of this theory in some detail and connect it to our work from previous sections.

Let $\mathcal{H}$ be a graded, connected Hopf algebra over $\mathbb{Q}$, and for $j = 0,1,\ldots, m-1$, let $\varphi^{(j)} : \mathcal{H} \to \mathbb{Q}$ be a character of $\mathcal{H}$. Let $\dot{\varphi} = (\varphi^{(0)}, \varphi^{(1)}, \ldots, \varphi^{(m-1)})$. Then the pair $(\mathcal{H}, \dot{\varphi})$ is what we will call an \emph{$m$-colored combinatorial Hopf algebra}. That is, we have our algebra $\mathcal{H}$ along with an (ordered) $m$-tuple of multiplicative linear functionals onto the base field of our algebra. In the category of $m$-colored combinatorial Hopf algebras, morphisms are Hopf homomorphisms $\psi: \mathcal{H} \to \mathcal{H}'$ such that $\varphi^{(j)} = (\varphi^{(j)})' \circ \psi$ for every $j = 0,1,\ldots, m-1$.

\begin{remark}\label{rem:err}
At first, this definition seems different from the one presented in \cite[Section~5.3]{BergeronHohlweg}, where a character is an algebra morphism $\varphi: \mathcal{H} \to \mathbb{Q}[\mathcal{C}_m]$. However, if we understand $\mathbb{Q}[\mathcal{C}_m]$ not as the usual group algebra, but as an algebra with multiplication given by \[\omega^j \omega^k = \begin{cases} \omega^j & \mbox{if } j = k,\\ 0 & \mbox{otherwise,} \end{cases}\] then their definition of character is naturally associated with ours.
\end{remark}

To put $\Q^{(m)}$ in this framework, let $\zeta^{(j)}_{\mathcal{Q}}: \Q^{(m)} \to \mathbb{Q}$ be the map that sends $x_{1,j}$ to 1 and all other variables to zero, extended algebraically. Then the pair $(\Q^{(m)},\dot{\zeta}_{\mathcal{Q}})$ is a colored combinatorial Hopf algebra. For each $j$ we have \[\zeta^{(j)}_{\mathcal{Q}}(M^{(m)}_\alpha) = \zeta^{(j)}_{\mathcal{Q}}(F^{(m)}_\alpha)=
\begin{cases}
1 & \text{if } \alpha = () \text{ or } \alpha = (\omega^j n) \\
0 & \text{otherwise.}
\end{cases}
\]
It is also interesting to consider the following character.
Let $\zeta_{\mathcal{Q}} = \zeta^{(0)}_{\mathcal{Q}}\zeta^{(1)}_{\mathcal{Q}}\cdots \zeta^{(m-1)}_{\mathcal{Q}}$, the convolution product of the $m$ characters taken in order. For the monomial and fundamental bases this gives
\[\zeta_{\mathcal{Q}}(M^{(m)}_\alpha) = \zeta_{\mathcal{Q}}(F^{(m)}_\alpha)=
\begin{cases}
1 & \text{if } \alpha = () \text{ or } \alpha = (\omega^{j_1} \alpha_1, \ldots, \omega^{j_k}\alpha_k), j_1 < \cdots < j_k, \\
0 & \text{otherwise.}
\end{cases}
\]

Before we present the universality theorem for colored combinatorial Hopf algebras, we first introduce the dual algebra of $\Q^{(m)}$, the colored noncommutative symmetric functions $\mathcal{N}Sym^{(m)} := \bigoplus_{n\geq 0} \mathcal{N}Sym^{(m)}_n$. See \cite{HivertNovelliThibon} or \cite{NovelliThibon} for more. This algebra is freely generated by the elements $H_{n}^{(j)}$, over all $n \geq 0$ and all $j = 0,1,\ldots,m-1$. A basis for $\mathcal{N}Sym^{(m)}_n$ is given by the set $\{ H_{\alpha} : \alpha \models_m n \}$, where if $\alpha = (\omega^{j_1} \alpha_1, \omega^{j_2}\alpha_2, \ldots, \omega^{j_k} \alpha_k)$ is an $m$-colored composition, then $H_{\alpha} := H_{\alpha_1}^{(j_1)} H_{\alpha_2}^{(j_2)} \cdots H_{\alpha_k}^{(j_k)}$. This $H$-basis is dual to the monomial basis of $\Q^{(m)}_n$, so that $(H_{\alpha})^* \leftrightarrow M_{\alpha}$.

The pair $(\Q^{(m)}, \dot{\zeta}_{\mathcal{Q}})$ satisfies the following universal property, analogous to the 1-colored case. See \cite[Theorem~5.4]{BergeronHohlweg} for a similar claim.

\begin{thm}\label{thm:colorUniv}
For any colored combinatorial Hopf algebra $(\mathcal{H}, \dot{\varphi})$, there is a unique morphism of colored combinatorial Hopf algebras, $\Psi: (\mathcal{H},\dot{\varphi}) \to (\Q^{(m)},\dot{\zeta}_{\mathcal{Q}})$. This map is given explicitly by
\begin{equation}\label{eq:explicit}
\Psi(h) = \sum_{ \alpha \models_m n} \varphi_{\alpha}(h)M^{(m)}_{\alpha},
\end{equation}
for any $h \in \mathcal{H}_n$, where for $\alpha = (\omega^{j_1}\alpha_1, \ldots, \omega^{j_k}\alpha_k)$, $\varphi_{\alpha}$ is the composite map \[\mathcal{H} \xrightarrow{\Delta^{(k-1)}} \mathcal{H}^{\otimes k} \twoheadrightarrow \mathcal{H}_{\alpha_1} \otimes \cdots \otimes \mathcal{H}_{\alpha_k} \xrightarrow{\textsf{m} \circ (\varphi^{(j_1)} \otimes \cdots \otimes \varphi^{(j_k)})} \mathbb{Q}.\]
\end{thm}

\begin{proof}
This proof is a straightforward generalization of the argument of \cite[Theorem~4.1]{AguiarBergeronSottile}. First consider $\mathcal{H}$ as a coalgebra. For every $n$ and $j$, let \[ \varphi_n^{(j)} := \varphi^{(j)}|_{\mathcal{H}_n} : \mathcal{H}_n \to \mathbb{Q}.\] Then $\varphi_n^{(j)} \in \mathcal{H}_n^*$, where $\mathcal{H}^* := \bigoplus_{n \geq 0} \mathcal{H}_n^*$ is the graded dual of $\mathcal{H}$. Since $\mathcal{N}Sym^{(m)}$ is freely generated by the $H_n^{(j)}$, there is a unique algebra homomorphism $\Phi: \mathcal{N}Sym^{(m)} \to \mathcal{H}^*$ such that \[ \Phi(H_n^{(j)}) = \varphi^{(j)}_n \] for every $n \geq 1$, $j = 0,1,\ldots, m-1$.

Let $\Psi := \Phi^*: \mathcal{H} \to \Q^{(m)}$ be the dual morphism of coalgebras, unique by the uniqueness of $\Phi$. For any $j$, we see that \[ (\zeta^{(j)}_{\mathcal{Q}} \circ \Psi)|_{\mathcal{H}_n} = \zeta^{(j)}_{\mathcal{Q}}|_{\Q^{(m)}_n} \circ \Psi|_{\mathcal{H}_n} = H_n^{(j)} \circ \Psi|_{\mathcal{H}_n} = \Phi(H_n^{(j)}) = \varphi_n^{(j)} = \varphi^{(j)}|_{\mathcal{H}_n}\] for every $n$. Therefore $\zeta^{(j)}_{\mathcal{Q}} \circ \Psi = \varphi^{(j)}$ for every $j$. We have showed that $\Psi$ is a morphism of colored combinatorial coalgebras. (Indeed this proves universality of $\Q^{(m)}$ in this category.)

We now check that $\Psi$ is a morphism of algebras, i.e., that $\Psi(gh) = \Psi(g)\Psi(h)$ for any $g,h \in \mathcal{H}$. Notice that $(\mathcal{H}\otimes\mathcal{H}, (\varphi^{(0)}\otimes\varphi^{(0)},\ldots, \varphi^{(m-1)}\otimes \varphi^{(m-1)}))$ is itself an $m$-colored combinatorial Hopf algebra. Now consider the composite maps \[\mathcal{H}\otimes\mathcal{H} \xrightarrow{\textsf{m}} \mathcal{H} \xrightarrow{\Psi} \Q^{(m)}\] and \[\mathcal{H}\otimes\mathcal{H} \xrightarrow{\Psi\otimes \Psi} \Q^{(m)}\otimes\Q^{(m)} \xrightarrow{\textsf{m}} \Q^{(m)},\] where \textsf{m} denotes the appropriate multiplication map. It is straightforward to check both $\Psi \circ \textsf{m}$ and $\textsf{m} \circ \Psi \otimes \Psi$ are morphisms of colored combinatorial coalgebras. But the previous paragraph shows that such a morphism to $\Q^{(m)}$ must be unique. Hence, $\Psi \circ \textsf{m}= \textsf{m} \circ \Psi \otimes \Psi$ and $\Psi$ is a morphism of combinatorial Hopf algebras. We have proved universality, now it remains to check our formula \eqref{eq:explicit}.

Since $\Phi(H_n^{(j)}) = \varphi_n^{(j)}$ and $\Phi$ is an algebra homomorphism, we have \[\Phi(H_{\alpha}) = \varphi_{\alpha_1}^{(j_1)}\varphi_{\alpha_2}^{(j_2)} \cdots \varphi_{\alpha_k}^{(j_k)} = \varphi_{\alpha}.\] Taking the dual, we get \[\Psi(h) = \sum_{ \alpha \models_m n} \varphi_{\alpha}(h)M^{(m)}_{\alpha},\] as desired.
\end{proof}

For a colored combinatorial Hopf algebra $(\mathcal{H},\dot{\varphi})$, let $\varphi = \varphi^{(0)} \varphi^{(1)} \cdots \varphi^{(m-1)}$ denote the convolution product of the characters of $\mathcal{H}$ (note that order matters). Then Theorem \ref{thm:colorUniv} implies that we have \[ \varphi = \zeta_{\mathcal{Q}} \circ \Psi.\] This fact leads to some nice formulas.

An $m$-colored poset $P$ will be called \emph{$j$-monochromatic} if all its labels share the same color, $\omega^j$. A colored poset is \emph{naturally labeled} if there exists a linear extension $\pi\in\mathcal{L}(P)$ such that its colored descent composition $C(\pi) = (\omega^{j_1}\alpha_1, \ldots, \omega^{j_k} \alpha_k)$ has the property that $j_1 < \cdots < j_k$, i.e., $\pi$ has no descents and its colors weakly increase. Notice that if there exists such a linear extension, it is unique. For $j=0,1,\ldots,m-1$, define the linear map $\zeta^{(j)}_{\mathcal{P}}:\mathcal{P}^{(m)}\to\mathbb{Q}$ by
\[\zeta^{(j)}_{\mathcal{P}}(P) = \begin{cases}
1 & \text{if $P$ is $j$-monochromatic and naturally labeled,}\\
0 & \text{otherwise}.
\end{cases}\]
It is easy to check that $\zeta^{(j)}_{\mathcal{P}}$ is a character, and that $(\mathcal{P}^{(m)}, \dot{\zeta_{\mathcal{P}}})$ is a colored combinatorial Hopf algebra.

If $\zeta_{\mathcal{P}} =\zeta^{(0)}_{\mathcal{P}} \zeta^{(1)}_{\mathcal{P}} \cdots \zeta^{(m-1)}_{\mathcal{P}}$, then we get
\[\zeta_{\mathcal{P}}(P) = \begin{cases}
1 & \text{if $P$ is naturally labeled,}\\
0 & \text{otherwise}.
\end{cases}\]
Moreover, an immediate consequence of \eqref{eq:mFTPP} and the definition of $\zeta^{(j)}_{\mathcal{P}}$ is that
\[\zeta^{(j)}_{\mathcal{P}} = \zeta^{(j)}_{\mathcal{Q}} \circ \Gamma_m\] for any $j = 0,1,\ldots,m-1$ and hence \[\zeta_{\mathcal{P}} = \zeta_{\mathcal{Q}} \circ \Gamma_m. \]
Thus $\Gamma_m$ is the unique homomorphism induced by $\dot{\zeta_{\mathcal{P}}}$.

We now try to understand $\Lambda_m$ from the perspective of $m$-colored combinatorial Hopf algebras. As a colored combinatorial Hopf algebra, the algebra of $m$-colored peak functions is the pair $(\mathbf{\Pi}^{(m)}, \dot{\zeta}_{\mathcal{Q}})$. Let $\pi \in \mathcal{C}_m \wr \mathfrak{S}_n$ and let $\widehat{C}(\pi) = (\omega^{j_1} \alpha_1, \ldots, \omega^{j_k} \alpha_k)$. Considering the number of colored enriched $\pi$-partitions $f$ of the form
\begin{align*}
(|f(\pi_1)|, |f(\pi_2)|, & \ldots, |f(\pi_n)|) \\
 & = ( \underbrace{\omega^{j_1} 1,\omega^{j_1} 1,\ldots, \omega^{j_1} 1}_{\alpha_1}, \underbrace{ \omega^{j_2} 1,\omega^{j_2} 1,\ldots , \omega^{j_2} 1}_{\alpha_2}, \ldots, \underbrace{\omega^{j_l} 1,\omega^{j_l} 1,\ldots, \omega^{j_l} 1}_{\alpha_k} ),
\end{align*}
we have
\begin{equation}\label{eq:mKzeta}
\zeta^{(j)}_{\mathcal{Q}}(K^{(m)}_{\alpha}) =
\begin{cases}
1 & \mbox{if } \alpha = (), \\
2 & \mbox{if } \alpha = (\omega^{j} n),\\
0 & \mbox{otherwise.}
\end{cases}
\end{equation}
and
\[
\zeta_{\mathcal{Q}}(K^{(m)}_{\alpha}) =
\begin{cases}
1 & \mbox{if } \alpha = (), \\
2^k & \mbox{if } \alpha = (\omega^{j_1} \alpha_1, \ldots, \omega^{j_k}\alpha_k), j_1 < \cdots < j_k,\\
0 & \mbox{otherwise.}
\end{cases}
\]

The {\em odd subalgebra} of an $m$-colored combinatorial Hopf algebra $(\mathcal{H},\dot{\varphi})$ is defined, as in \cite[Definition~5.1]{AguiarBergeronSottile}, to be the largest graded subcoalgebra $S_-(\mathcal{H},\dot{\varphi})$ of $\mathcal{H}$ such that for all $h\in\mathcal{H}$, $\overline{\varphi}^{(j)}(h) = (\varphi^{(j)})^{-1}(h)$ for $j=0,1,\ldots, m-1.$ It is straightforward to adapt the techniques from \cite[Section~6]{AguiarBergeronSottile} to show that $\mathbf{\Pi}^{(m)} = S_-(\Q^{(m)},\dot{\zeta}_{\mathcal{Q}})$ and $(\mathbf{\Pi}^{(m)},\dot{\zeta}_{\mathcal{Q}})$ is the terminal object for odd $m$-colored combinatorial Hopf algebras. Moreover, $\Theta_m$ is the unique Hopf homomorphism induced by the odd characters $\nu^{(j)}_{\mathcal{Q}} := (\bar{\zeta}^{(j)}_{\mathcal{Q}})^{-1}  \zeta^{(j)}_{\mathcal{Q}}$, $j=0,1,\ldots, m-1$. In other words, $\nu^{(j)}_{\mathcal{Q}} = \zeta^{(j)}_{\mathcal{Q}} \circ \Theta_m$ for every $j$ and \[\nu_{\mathcal{Q}} := \nu^{(0)}_{\mathcal{Q}}\nu^{(1)}_{\mathcal{Q}} \cdots \nu^{(m-1)}_{\mathcal{Q}} = \zeta_{\mathcal{Q}} \circ \Theta_m.\]
We leave the details to the reader.

For each $j$, consider the odd character $\nu^{(j)}_{\mathcal{P}} =  (\bar{\zeta}^{(j)}_{\mathcal{P}})^{-1}  \zeta^{(j)}_{\mathcal{P}}$. Since $\zeta^{(j)}_{\mathcal{P}} = \zeta^{(j)}_{\mathcal{Q}} \circ \Gamma_m$ and $\Gamma_m$ is a Hopf homomorphism, it follows that $\nu^{(j)}_{\mathcal{P}} = \nu^{(j)}_{\mathcal{Q}} \circ \Gamma_m$, and so
\begin{equation}\label{eq:mnu_P}
 \nu^{(j)}_{\mathcal{P}} = \zeta^{(j)}_{\mathcal{Q}} \circ \Lambda_m,
 \end{equation}
which also implies \[\nu_{\mathcal{P}} :=\nu^{(0)}_{\mathcal{P}}\nu^{(1)}_{\mathcal{P}} \cdots \nu^{(m-1)}_{\mathcal{P}} = \zeta_{\mathcal{Q}} \circ \Lambda_m.\]
Therefore $\Lambda_m$ is the unique homomorphism induced by $\dot{\nu}_{\mathcal{P}}$.

Combining \eqref{eq:mnu_P}, \eqref{eq:mKzeta}, and \eqref{eq:mLambdasum}, we see that $\nu^{(j)}_{\mathcal{P}}(P)$ and $\nu_{\mathcal{P}}(P)$ depend on those linear extensions with weakly increasing color and no peaks.
\begin{prp}
For any non-empty labeled poset $P$,
\[ \nu^{(j)}_{\mathcal{P}}(P) = 2 \cdot \#\{\pi\in\mathcal{L}(P)~|~\Pe(\pi) = \emptyset \mbox{ and } \varepsilon(\pi_i) = j \mbox{ for all } i \}, \]
and
\[ \nu_{\mathcal{P}}(P) = 2^k \cdot \#\{\pi\in\mathcal{L}(P)~|~\Pe(\pi) = \emptyset \mbox{ and } \varepsilon(\pi_i) \leq \varepsilon(\pi_{i+1}) \mbox{ for all } i\},\] where $k$ is the number of distinct colors of $P$.
\end{prp}


\begin{thebibliography}{99}


\bibitem{AguiarBergeronSottile} M. Aguiar, N. Bergeron, and F. Sottile, \emph{Combinatorial Hopf algebras and generalized Dehn-Sommerville relations}, Compositio Mathematica {\bf 142} (2006), 1--30.

\bibitem{AguiarSottile}  M. Aguiar and F. Sottile, \emph{Structure of the Malvenuto-Reutenauer Hopf algebra of permutations},  Advances in Mathematics  {\bf 191}  (2005), 225--275.

\bibitem{BaumannHohlweg} P. Baumann and C. Hohlweg, \emph{A Solomon descent theory for the wreath products $G\wr S_n$}, Transaction of the AMS, to appear.

\bibitem{BergeronHohlweg} N. Bergeron and C. Hohlweg, \emph{Coloured peak algebras and Hopf algebras}, Journal of Algebraic Combinatorics {\bf 24}  (2006), 299--330.

\bibitem{BergeronMSW} N. Bergeron, S. Mykytiuk, F. Sottile, and S. van Willigenburg, \emph{Shifted quasi-symmetric functions and the Hopf algebra of peak functions}, Discrete Mathematics {\bf 256} (2002), 57--66.

\bibitem{Chow} C.-O. Chow, \emph{Noncommutative symmetric functions of type B},
Ph.D. thesis, MIT (2001).

\bibitem{Ehrenborg} R. Ehrenborg, \emph{On posets and Hopf algebras}, Advances in Mathematics {\bf 119} (1996), 1--25.

\bibitem{Gessel} I. Gessel, \emph{Multipartite $P$-partitions and
inner products of skew Schur functions}, Contemporary Mathematics
{\bf 34} (1984), 289--317.

\bibitem{HivertNovelliThibon} F. Hivert, J.-C. Novelli, and J.-Y. Thibon, \emph{Yang-Baxter bases of $0$-Hecke algebras and representation theory of $0$-Ariki-Koike-Shoji algebras}, math.CO/0506546v1.


\bibitem{HsiaoPetersen} S. Hsiao and T. K. Petersen, \emph{The Hopf algebras of type B quasisymmetric functions and peak functions}, math.CO/0610976.

\bibitem{JoniRota} S. Joni and G.-C. Rota, \emph{Coalgebras and bialgebras in combinatorics}, Stud. Appl. Math. {\bf 61} (1979), 93--139.

\bibitem{Malvenuto} C. Malvenuto, \emph{Produits et coproduits des fonctions quasi-sym\'etriques et de l'alg\`ebre des descents}, Ph.D. Thesis, Universit\'e du Qu\'ebec \`a Montr\'eal, 1993.

\bibitem{MalvenutoReutenauer} C. Malvenuto and C. Reutenauer, \emph{Duality between quasi-symmetric functions and the Solomon descent algebra}, Journal of Algebra {\bf 177} (1995), 967--982.

\bibitem{MilnorMoore} J. Milnor and J. Moore, \emph{On the structure of Hopf Algebras}, Annals of Mathematics {\bf 81} (1965) 211--264.

\bibitem{NovelliThibon} J.-C. Novelli and J.-Y. Thibon, \emph{Free quasi-symmetric functions of arbitrary level}, math.CO/0405597.

\bibitem{Nyman} K. Nyman, \emph{The peak algebra of the symmetric group}, Journal of Algebraic Combinatorics {\bf 17} (2003), 309--322.

\bibitem{Petersen1} T. K. Petersen, \emph{A note on three types of quasisymmetric functions},  Electronic Journal of Combinatorics {\bf 12}  (2005), R61.

\bibitem{Petersen2} T. K. Petersen, \emph{Enriched $P$-partitions and peak algebras}, Advances in Mathematics, to appear.

\bibitem{Poirier} S. Poirier, \emph{Cycle type and descent set in wreath products}, Discrete Mathematics {\bf 180} (1998), 315--343.



\bibitem{Rota} G.-C. Rota, \emph{On the foundations of combinatorial theory. I. Theory of M\"{o}bius functions}, Z. Wahrscheinlichkeitstheorie Verw. Gebiete {\bf 2} (1964), 340--368.

%

\bibitem{Solomon} L. Solomon, \emph{A Mackey formula in the group ring of a finite Coxeter group}, Journal of Algebra {\bf 41} (1976), 255--264.

\bibitem{Stanley1} R. Stanley, \emph{Enumerative Combinatorics, Volume I}, Cambridge University Press, 1997.

\bibitem{Stanley2} R. Stanley, \emph{Enumerative Combinatorics, Volume II}, Cambridge University Press, 2001.

\bibitem{Stembridge} J. Stembridge, \emph{Enriched ${P}$-partitions}, Transactions of the American Mathematical Society {\bf 349} (1997), 763--788.

\end{thebibliography}
\end{document}